\documentclass{amsart}
\usepackage{amssymb} 
\usepackage{color}
\usepackage[all]{xy}
\usepackage{enumerate}

 \newtheorem{theorem}{Theorem}[section]

 \newtheorem{proposition}[theorem]{Proposition}
 \newtheorem{lemma}[theorem]{Lemma}
 \newtheorem{corollary}[theorem]{Corollary}
 \newtheorem{remark}[theorem]{Remark}
 \newtheorem{definition}[theorem]{Definition}
 
 \newtheorem{example}[theorem]{Example}

 \newtheorem{condition}[theorem]{Condition}
 
 \numberwithin{equation}{section}

 \def\subrel#1#2{\mathrel{\mathop{#2}\limits_{#1}}}

\def\cat{{\rm Cat}_{\infty}}

\def\overcat#1{\cat\mbox{$\scriptstyle /#1$}}
\def\duo{{\rm Duo}_{\infty}}
\def\wcat{\widehat{\rm Cat}_{\infty}}

\def\nonoperad{{\rm Op}_{\infty}^{\rm ns}}
\def\nonoperad{{\rm Op}_{\infty}^{\rm ns}}

\def\monno#1{{\rm Mon}(#1)}

\def\laxmoncat{\mathsf{Mon}^{\rm lax}(\cat)}

\def\laxmon{\mathsf{Mon}^{\rm lax}}
\def\oplaxmon{\mathsf{Mon}^{\rm oplax}}
\def\oplaxmoncat{\mathsf{Mon}^{\rm oplax}(\cat)}

\def\coalg#1{{}_{#1}{\rm cAlg}_{#1}}
\def\alg{{\rm Alg}}
\def\bialg#1{{}_{#1}{\rm bMon}_{#1}}
\def\dalg#1{{\rm dMon}_{#1}}
\def\dcalg#1{{\rm dcMon}_{#1}}
\def\cart#1{{\rm Cart}\mbox{$\scriptstyle /#1$}}
\def\cocart#1{{\rm coCart}\mbox{$\scriptstyle /#1$}}
\def\mixed#1{{\rm Mfib}\mbox{$\scriptstyle /#1$}}

\begin{document}

\title%[On Duoidal $\infty$-categories]
{On duoidal $\infty$-categories}
%{Duoidal $\infty$-categories, I:\\
%The $\infty$-category of duoidal $\infty$-categories\\ 
%and bilax monoidal functors}

\author{Takeshi Torii}
\address{Department of Mathematics, 
%Faculty of Sience, 
Okayama University,
Okayama 700--8530, Japan}
\email{torii@math.okayama-u.ac.jp}
\thanks{The author was partially supported by
JSPS KAKENHI Grant Numbers JP17K05253 and JP23K03113.}

%\subjclass[2000]{18D10 (primary), 18D35, 18D50, 55U40 (secondary)}

\date{January 25, 2025 (version~5.3)}

\begin{abstract}
A duoidal category is a category equipped with
two monoidal structures in which one is (op)lax
monoidal with respect to the other.
In this paper
we introduce duoidal $\infty$-categories
which are counterparts of duoidal categories 
in the setting of $\infty$-categories.
There are three kinds of functors
between duoidal $\infty$-categories,
which are called bilax, double lax, and double oplax monoidal functors.
We make three formulations of 
$\infty$-categories of duoidal $\infty$-categories
according to which functors we take.
Furthermore, 
corresponding to the three kinds 
of functors,
we define bimonoids, double monoids,
and double comonoids in duoidal
$\infty$-categories.
\end{abstract}

\maketitle

\section{Introduction}

The goal of this paper is to introduce duoidal categories
in the setting of $\infty$-categories.
A duoidal category is a category equipped with
two monoidal products in which they are compatible 
in the following sense:
the functors and natural transformations
defining one monoidal structure 
are (op)lax monoidal with respect to the other monoidal structure.

The notion of duoidal categories
have been introduced in \cite{Aguiar-Mahajan}
under the name of $2$-monoidal categories.
Aguiar and Mahajan~\cite{Aguiar-Mahajan} 
have developed a basic theory of duoidal categories
and three kinds of functors between them.
The name of duoidal categories has appeared
in \cite{Street}.
Inspired by $3$-manifold invariants,
Street~\cite{Street} has developed 
a theory of duoidal categories
and centers of monoids in monoidal categories.
After that,
duoidal categories have been used in many papers.
In \cite{Batanin-Markl} Batanin and Markl
have considered a theory of centers and homotopy centers
of monoids in monoidal categories which are 
enriched in duoidal categories,
and in \cite{Batanin-Markl2} they 
have proved duoidal Deligne conjecture using 
their previous work.
Booker and Street have further developed
a general theory of duoidal categories
in \cite{Booker-Street}.
We can consider bimonoids in a duoidal category,
which is a generalization of bialgebras
in a symmetric monoidal category.
In order to generalize the notion of Hopf algebras,
B\"{o}hm, Chen, and Zhang
have studied Hopf monoids and the fundamental theorem of Hopf modules
in duoidal categories
in \cite{BCZ}.
Furthermore,
B\"{o}hm and Lack 
have studied a notion of antipode on bimonoids
in duoidal categories under some circumstances
in \cite{Bohm-Lack}.

There are many approaches to the theory of
higher categories, and
we have many models of $(\infty,1)$-categories,
where $(\infty,1)$-categories are $\infty$-categories
in which all $k$-morphisms are invertible for $k>1$.
Quasi-categories are one of these models,
which were originally
introduced by Boardman-Vogt~\cite{Boardman-Vogt}
in their study of homotopy invariant
algebraic structures.
After that,
quasi-categories are extensively developed
by Joyal~\cite{Joyal1,Joyal2},
Lurie~\cite{Lurie1,Lurie2},
and many others.
Quasi-categories and other models
of $(\infty,1)$-categories have many applications in many areas
of mathematics.
According to Lurie~\cite{Lurie1},
we use the term $\infty$-category
to refer to quasi-categories and
other models of $(\infty,1)$-categories.

Comodules over Hopf algebroids of co-operations
for generalized (co)homology theories 
play fundamental roles in chromatic homotopy theory~
\cite{MRW, Morava, Ravenel}.
We can regard Hopf algebroids as
examples of bimonoids
in the duoidal category of bimodules
over graded commutative rings. 
In \cite{Torii}
we have studied $\infty$-categories of comodule spectra
over coalgebras of co-operations for generalized (co)homology theories. 
Although we did not use multiplicative structures
on $\infty$-categories of comodule spectra
in \cite{Torii},
we realized that
we need to develop a theory
of duoidal categories in the setting of $\infty$-categories
in order to discuss multiplicative structures
on these $\infty$-categories of comodule spectra. 

In this paper we will introduce duoidal categories
in the setting of $\infty$-categories,
which we call duoidal $\infty$-categories.
We have three kinds of functors between 
duoidal categories,
which are called 
bilax, double lax, and double oplax monoidal functors.
Corresponding
to theses functors,
we define three kinds of functors
between duoidal $\infty$-categories.
We construct three $\infty$-categories
\[ \duo^{\rm bilax},\quad \duo^{\rm dlax},\quad
   \duo^{\rm doplax}\]
of duoidal $\infty$-categories with
bilax, double lax, and double oplax monoidal functors,
respectively.

There are notions of 
bimonoids, double monoids, and double comonoids 
in duoidal categories.
They are generalizations
of bialgebras, commutative algebras,
cocommutative coalgebras in braided monoidal categories,
respectively.
We will introduce counterparts 
in duoidal $\infty$-categories, 
and construct three $\infty$-categories
\[ \bialg{}(X), \quad \dalg{}(X),\quad \dcalg{}(X) \]
of bimonoids, double monoids, and double comonoids
in a duoidal $\infty$-category $X$.
These constructions determine three functors
\[ \begin{array}{rlcl}
   \bialg{}:& \duo^{\rm bilax}&\longrightarrow& \cat,\\[2mm]
   \dalg{}:& \duo^{\rm dlax}&\longrightarrow& \cat, \\[2mm]
   \dcalg{}:&\duo^{\rm doplax}&\longrightarrow& \cat.\\
   \end{array}\]

For a duoidal $\infty$-category $X$
with two monoidal structures $\boxtimes$
and $\otimes$,
where $\otimes$ is lax monoidal
with respect to $\boxtimes$,
we will show that the $\infty$-category
$\coalg{}{}^{\otimes}(X)$
of coalgebra objects of the monoidal $\infty$-category
$(X,\otimes)$ 
supports a monoidal structure
induced by $\boxtimes$
(Theorem~\ref{thm:cP-monoidal-infinity-category}).
This construction determines a functor
$\coalg{}{}^{\otimes}: \duo^{\rm bilax}\to\laxmoncat$,
where $\laxmoncat$ is the $\infty$-category
of monoidal $\infty$-categories and
lax monoidal functors.
Similarly,
we will show that
there is a functor
$\alg{}{}^{\boxtimes}: \duo^{\rm bilax}\to\oplaxmoncat$,
which associates to a duoidal $\infty$-category
$X$ the $\infty$-category of algebra
objects of the monoidal $\infty$-category $(X,\boxtimes)$,
where $\oplaxmoncat$
is the $\infty$-category of monoidal
$\infty$-categories and oplax monoidal functors
(Remark~\ref{remark:algoduo-bilax-oplax-monoidal-functor}). 

One of our main theorems is as follows. 

\begin{theorem}[Theorem~\ref{theorem:bimonoid-is-algebra-of-coalgebra},
Corollary~\ref{cor:equivalence-bimod-alg-coalg},
and Remark~\ref{remark:coalg-alg-dual-duoidal}]
\label{theorem:one-of-main-theorems}
For a duoidal $\infty$-category $X$,
there are natural equivalences
\[ \begin{array}{rcl}
    \bialg{}(X)&\simeq& \alg{}(\coalg{}{}^{\otimes}(X))\\[2mm]
               &\simeq& \coalg{}(\alg{}{}^{\boxtimes}(X))\\
   \end{array}  \]
of $\infty$-categories.
As a result,
there are equivalences
\[ \begin{array}{rcl}
     \bialg{}&\simeq& \alg{}\circ\coalg{}^{\otimes}\\[2mm]
             &\simeq& \coalg{}\circ\alg{}{}^{\boxtimes}\\
   \end{array}  \]
of functors from
$\duo^{\rm bilax}$ to $\cat$.
\end{theorem}

We will prove similar results on 
the functors $\dalg{}: \duo^{\rm dlax}\to\cat$
and $\dcalg{}:\duo^{\rm doplax}\to\cat$
(Theorems~\ref{theorem:duoidal-alg-monoidal-double-formulation},
\ref{theorem:dalg-equivalent-alg-alg-boxtimes},
and \ref{theorem:dcoalg-equivalent-coalg-coalg-boxtimes}).

Based on the results in this paper,
we will generalize duoidal $\infty$-categories
to higher monoidal $\infty$-categories in \cite{Torii2}.
In \cite{Torii4}
we will give an example of duoidal $\infty$-categories
of operadic modules.
In \cite{Torii5}
we will give another example of duoidal $\infty$-categories
obtained from map monoidales in monoidal
$\infty$-bicategories.

The organization of this paper
is as follows:
%{\color{red}
In \S\ref{section:review-duoidal}
we recall the notion of duoidal categories
and three kinds of functors between them
in the classical setting.
%}
In \S\ref{section:monoidal-infinity-categories}
we review monoidal $\infty$-categories and functors between them.
We recall lax and oplax monoidal functors
between monoidal $\infty$-categories.
In \S\ref{section:mixed-fibrations}
we introduce mixed fibrations and study their properties.
In \S\ref{subsection:Functor-infinity-categories}
we study monoid objects of the $\infty$-category
of a slice $\infty$-category $\overcat{T}$.
We consider lax and oplax monoidal functors
between monoid objects of $\overcat{T}$.
Motivated by 
the results in \S\ref{subsection:Functor-infinity-categories},
we introduce mixed fibrations
in \S\ref{subsection:mixed-fibration}.
Although mixed fibrations
have an obvious duality by definition,
we also give another asymmetric description 
of mixed fibrations.
%In \S\ref{subsection:marked-mixed-fibrations}
We also study mixed fibrations over marked simplicial sets.
In \S\ref{section:duoidal-infty-categories}
we introduce duoidal $\infty$-categories
as mixed fibrations which satisfy the Segal conditions.
%{\color{red}
As in the case of classical setting,
there are three kinds of functors between duoidal 
$\infty$-categories.
%}
According to which kind of functors of duoidal $\infty$-categories
we take,
we will give three formulations
of $\infty$-categories of duoidal $\infty$-categories.
In \S\ref{section:bialgebra}
we introduce notions of 
bimonoids, double monoids, and double comonoids
in duoidal $\infty$-categories.
In \S\ref{subsection:algebra-coalgebra-monoida-category}
we recall algebra and coalgebra objects
in monoidal $\infty$-categories.
In \S\ref{subsection:bimonoida-in-duoidal-infty-category}
we define bimonoids in duoidal $\infty$-categories
and prove Theorem~\ref{theorem:one-of-main-theorems}.
In \S\ref{subsection:double-monoid-comonoid-in-duoidal}
we define double monoids and
double comonoids in duoidal $\infty$-categories
and prove similar results on them
to Theorem~\ref{theorem:one-of-main-theorems}.

\begin{remark}\rm
In this paper we think ordinary categories
are a special kind of $\infty$-categories.
Thus, we do not distinguish notationally
between ordinary categories and their nerves.
In particular,
we identify the simplicial indexing category $\Delta$
with its nerve $N(\Delta)$.
\end{remark}

%{\color{red}
\noindent
{\bf Acknowledgements}.
The author would like to thank
the referee for his/her
valuable comments and useful suggestions.
%}

%{\color{red}
%\subsection{Review of duoidal categories}
\section{Duoidal categories in the classical setting}
\label{section:review-duoidal}
%}

The notion of duoidal categories 
has been introduced in \cite{Aguiar-Mahajan}
by the name of $2$-monoidal categories,
and the name of duoidal categories 
has appeared in \cite{Street}.
In this section we recall 
duoidal categories in the classical setting
and three kinds of functors between them.

\if0
First,
we recall the definition of monoidal category.

\begin{definition}
[{cf.~\cite[Definition~1.1]{Aguiar-Mahajan}}]
A monoidal category is a triple
$(\mathcal{C},\otimes,1_{\otimes})$,
where $\mathcal{C}$ is a category,
$\otimes$ is a functor
\[ \otimes: \mathcal{C}\times\mathcal{C}
\longrightarrow \mathcal{C},\]
\end{definition}
\fi

\begin{definition}
%[{cf.~\cite{Aguiar-Mahajan, Batanin-Markl, BCZ, 
%Bohm-Lack, Booker-Street, Street}}]\rm
%{\color{red}
[{cf.~\cite[Definition~6.1]{Aguiar-Mahajan}}]
%}
\rm
A duoidal category is a quintuple
$(\mathcal{D},\boxtimes,1_{\boxtimes},\otimes,1_{\otimes})$,
where
%\begin{enumerate}
%\item
$(\mathcal{D},\boxtimes,1_{\boxtimes})$
and
$(\mathcal{D},\otimes,1_{\otimes})$
are monoidal categories,
along with
%\item
a natural transformation
\begin{equation}\label{eq:interchange-law}
    \zeta_{A,B,C,D}:
    (A\otimes B)\boxtimes(C\otimes D)\longrightarrow 
    (A\boxtimes C)\otimes(B\boxtimes D),
\end{equation}
and
%\item
three morphisms
\begin{equation}\label{eq:three-morphisms}
  \Delta:1_{\boxtimes}\to  1_{\boxtimes}\otimes 1_{\boxtimes},\qquad
   \mu: 1_{\otimes}\boxtimes 1_{\otimes} \to 1_{\otimes},\qquad 
   \iota=\epsilon:1_{\boxtimes}\to 1_{\otimes}
\end{equation}
%\end{enumerate}
satisfying the following conditions\mbox{\rm :}

\noindent
{\bf Associativity.}
The following diagrams commute
\[ \xymatrix{
    ((A\otimes B)\boxtimes (C\otimes D))
    \boxtimes (E\otimes F)
    \ar[r]^{\alpha}\ar[d]_{\zeta\boxtimes{\rm id}} &
    (A\otimes B)\boxtimes
    ((C\otimes D)\boxtimes (E\otimes F))
    \ar[d]_{{\rm id}\boxtimes\zeta}\\
    ((A\boxtimes C)\otimes (B\boxtimes D))
    \boxtimes (E\otimes F)
    \ar[d]_{\zeta}
    & (A\otimes B)\boxtimes((C\boxtimes E)\otimes
    (D\boxtimes F))
    \ar[d]^{\zeta}\\
    ((A\boxtimes C)\boxtimes E)\otimes
    ((B\boxtimes D)\boxtimes F)
    \ar[r]_{\alpha\otimes\alpha}&
    (A\boxtimes (C\boxtimes E))\otimes
    (B\boxtimes (D\boxtimes F)),\\    
}\]

\[ \xymatrix{
  ((A\otimes B)\otimes C)\boxtimes
  ((D\otimes E)\otimes F)
  \ar[r]^{\alpha\boxtimes\alpha}\ar[d]^{\zeta}&
  (A\otimes (B\otimes C))\boxtimes
  (D\otimes(E\otimes F))
  \ar[d]^{\zeta}\\
  ((A\otimes B)\boxtimes (D\otimes E))\otimes
  (C\boxtimes F)\ar[d]^{\zeta\otimes{\rm id}} &
  (A\boxtimes D)\otimes
  ((B\otimes C)\boxtimes (E\otimes F))
  \ar[d]^{{\rm id}\otimes\zeta}\\
  ((A\boxtimes D)\otimes (B\boxtimes E))
  \otimes (C\boxtimes F)\ar[r]^{\alpha}&
  (A\boxtimes D)\otimes
  ((B\boxtimes E)\otimes (C\boxtimes F)),\\
}\]    
where $\alpha$
is the associator
of either monoidal category.

\noindent
{\bf Unitality.}
The following diagrams commute
\[ \xymatrix{
  1_{\boxtimes}\boxtimes (A\otimes B)
  \ar[r]^-{\Delta\boxtimes{\rm id}}&
  (1_{\boxtimes}\otimes 1_{\boxtimes})\boxtimes
  (A\otimes B)\ar[d]^-{\zeta}\\
  A\otimes B\ar[u]^-{\lambda}
  \ar[r]_-{\lambda\otimes\lambda}&
  (1_{\boxtimes}\boxtimes A)\otimes
  (1_{\boxtimes}\boxtimes B),\\
   }\qquad
   \xymatrix{
  (A\otimes B)\boxtimes 1_{\boxtimes}
  \ar[r]^-{{\rm id}\boxtimes\Delta}&
  (A\otimes B)\ar[d]^-{\zeta}\boxtimes
  (1_{\boxtimes}\otimes 1_{\boxtimes})\\
  A\otimes B\ar[u]^-{\rho}
  \ar[r]_-{\rho\otimes\rho}&
  (A\boxtimes 1_{\boxtimes})\otimes
  (B\boxtimes 1_{\boxtimes}),\\
  }\]
\[ \xymatrix{
  1_{\otimes}\otimes (A\boxtimes B)  &
  (1_{\otimes}\boxtimes 1_{\otimes})\otimes
  (A\boxtimes B)\ar[l]_-{\mu\otimes{\rm id}}\\
  A\boxtimes B\ar[u]^-{\lambda}
  \ar[r]_-{\lambda\boxtimes\lambda}&
  (1_{\otimes}\otimes A)\boxtimes
  (1_{\otimes}\otimes B)
  \ar[u]_-{\zeta},\\
   }\qquad
   \xymatrix{
  (A\boxtimes B)\otimes 1_{\otimes}  &
  (A\boxtimes B)\otimes
  (1_{\otimes}\boxtimes 1_{\otimes})
  \ar[l]_-{{\rm id}\otimes\mu}\\
  A\boxtimes B\ar[u]^-{\rho}
  \ar[r]_-{\rho\boxtimes\rho}&
  (A\otimes 1_{\otimes})\boxtimes
  (B\otimes 1_{\otimes})
  \ar[u]_-{\zeta},\\
  }\]
where $\lambda$ is the left and $\rho$
is the right unitor
of either monoidal category.
%where $\lambda$ and $\rho$ are
%the left and right unitors
%of either monoidal category,
%respectively.

\noindent
{\bf Compatibility of units.}
The units $1_{\boxtimes}$ and $1_{\otimes}$
are compatible in the following sense:
\begin{itemize}
\item
$(1_{\otimes},\mu,\iota)$
is a monoid in $(\mathcal{D},\boxtimes,1_{\boxtimes})$,
and
\item
$(1_{\boxtimes}, \Delta, \epsilon)$
is a comonoid in $(\mathcal{D},\otimes, 1_{\otimes})$.
\end{itemize}
\end{definition}

There are several equivalent descriptions
for the notion of duoidal category.
One of them is given in terms of pseudomonoids
in monoidal 2-categories
(see \cite[Appendix~C]{Aguiar-Mahajan}
for monoidal 2-categories and pseudomonoids).
%}
For a monoidal $2$-category $\mathcal{M}$,
we can consider pseudomonoids in $\mathcal{M}$. 
Furthermore, there are notions
of lax and oplax morphisms between them.
We denote by 
\[ {\rm Mon}^{\rm lax}(\mathcal{M}) \]
the $2$-category whose $0$-cells are pseudomonoids in $\mathcal{M}$,
$1$-cells are lax morphisms between pseudomonoids, and  
$2$-cells are morphisms between lax morphisms.
We also denote by
\[ {\rm Mon}^{\rm oplax}(\mathcal{M}) \]
the $2$-category whose $0$-cells are pseudomonoids in $\mathcal{M}$,
$1$-cells are oplax morphisms between pseudomonoids, and  
$2$-cells are morphisms between oplax morphisms.

Let 
${\rm Cat}$
be the $2$-category
of categories, functors, and natural transformations.
Notice  that ${\rm Cat}$ is a monoidal $2$-category
under Cartesian product. 
Applying ${\rm Mon}^{\rm lax}(-)$ and ${\rm Mon}^{\rm oplax}(-)$
to ${\rm Cat}$,
we obtain $2$-categories
${\rm Mon}^{\rm lax}({\rm Cat})$ and
${\rm Mon}^{\rm oplax}({\rm Cat})$, respectively.
Furthermore,
${\rm Mon}^{\rm lax}({\rm Cat})$ and
${\rm Mon}^{\rm oplax}({\rm Cat})$
are again monoidal $2$-categories under Cartesian product.

\begin{proposition}[{cf.~\cite[Proposition~6.73]{Aguiar-Mahajan}}]
\label{prop:duoidal-pseudomonoid-identified}
A duoidal category is identified with
a pseudomonoid in ${\rm Mon}^{\rm lax}({\rm Cat})$.
A duoidal category is also identified with
a pseudomonoid in ${\rm Mon}^{\rm oplax}({\rm Cat})$.
\end{proposition}

Next,
we recall three kinds of functors between duoidal categories.
Since ${\rm Mon}^{\rm oplax}({\rm Cat})$
%and ${\rm Mon}^{\rm oplax}({\rm Mon}^{\rm oplax}(\cat))$
is a monoidal $2$-category,
we can consider a $2$-category
\[ {\rm Mon}^{\rm lax}({\rm Mon}^{\rm oplax}({\rm Cat})).\]
%and
%\[ {\rm Mon}^{\rm oplax}({\rm Mon}^{\rm lax}({\rm Cat}))\]
By Proposition~\ref{prop:duoidal-pseudomonoid-identified},
$0$-cells
of ${\rm Mon}^{\rm lax}({\rm Mon}^{\rm oplax}({\rm Cat}))$
are identified with duoidal categories.

\begin{definition}\rm
A bilax monoidal functor 
is a $1$-cell of 
${\rm Mon}^{\rm lax}({\rm Mon}^{\rm oplax}({\rm Cat}))$.
\end{definition}

In other words,
a bilax monoidal functor
between duoidal categories
$(\mathcal{C},\boxtimes,1_{\boxtimes},\otimes,1_{\otimes})$
and
$(\mathcal{D},\boxtimes,1_{\boxtimes},\otimes,1_{\otimes})$
is a triple $(F,\varphi,\psi)$,
where 
\begin{enumerate}
\item
$F: \mathcal{C}\to \mathcal{D}$
is a functor,
\item
$(F,\varphi): (\mathcal{C},\boxtimes,1_{\boxtimes})\to
(\mathcal{D},\boxtimes,1_{\boxtimes})$
is lax monoidal, and
\item
$(F,\psi): (\mathcal{C},\otimes,1_{\otimes})\to
(\mathcal{D},\otimes,1_{\otimes})$ 
is oplax monoidal, 
\end{enumerate} 
satisfying appropriate compatibility conditions
(see \cite[Definition~6.50]{Aguiar-Mahajan} for details).

Note that 
there is an equivalence
\[ {\rm Mon}^{\rm lax}({\rm Mon}^{\rm oplax}({\rm Cat}))
   \simeq
  {\rm Mon}^{\rm oplax}({\rm Mon}^{\rm lax}({\rm Cat})) \]
of $2$-categories 
by \cite[Proposition~6.75]{Aguiar-Mahajan}.
Hence we can identify bilax monoidal functors
with $1$-cells of the $2$-category 
${\rm Mon}^{\rm oplax}({\rm Mon}^{\rm lax}({\rm Cat}))$.

We can also consider a $2$-category
\[ {\rm Mon}^{\rm lax}({\rm Mon}^{\rm lax}({\rm Cat})).\]
%and
%\[ {\rm Mon}^{\rm oplax}({\rm Mon}^{\rm oplax}({\rm Cat})).\]

\begin{definition}\rm
A double lax monoidal functor is a $1$-cell of 
${\rm Mon}^{\rm lax}({\rm Mon}^{\rm lax}({\rm Cat}))$.
%A double oplax functor is a $1$-cell of
%${\rm Mon}^{\rm oplax}({\rm Mon}^{\rm oplax}(\cat))$.
%We denote by
%\[ {\rm Duo}^{\rm dlax} \]
%the $2$-category ${\rm Mon}^{\rm lax}({\rm Mon}^{\rm lax}(\cat))$
%and call it the $2$-category of duoidal categories,
%double lax functors, and morphisms of double lax functors.
\end{definition}

As in the case of bilax monoidal functors,
a double lax monoidal functor
between %duoidal categories
$(\mathcal{C},\boxtimes,1_{\boxtimes},\otimes,1_{\otimes})$
and
$(\mathcal{D},\boxtimes,1_{\boxtimes},\otimes,1_{\otimes})$
is a triple $(F,\varphi,\psi)$,
where 
\begin{enumerate}
\item
$F: \mathcal{C}\to \mathcal{D}$
is a functor,
\item
$(F,\varphi): (\mathcal{C},\boxtimes,1_{\boxtimes})\to
(\mathcal{D},\boxtimes,1_{\boxtimes})$
is lax monoidal, and
\item
$(F,\psi): (\mathcal{C},\otimes,1_{\otimes})\to
(\mathcal{D},\otimes,1_{\otimes})$ 
is lax monoidal, 
\end{enumerate} 
satisfying appropriate compatibility conditions
(see \cite[Definition~6.54]{Aguiar-Mahajan} for details).
%(see \cite[Definition~6.54 and 6.55]{Aguiar-Mahajan} for details).

Similarly,
we can consider a $2$-category
\[ {\rm Mon}^{\rm oplax}({\rm Mon}^{\rm oplax}(\cat)).\]

\begin{definition}\rm
%A double lax functor is a $1$-cell of 
%${\rm Mon}^{\rm lax}({\rm Mon}^{\rm lax}(\cat))$.
A double oplax monoidal functor is a $1$-cell of
${\rm Mon}^{\rm oplax}({\rm Mon}^{\rm oplax}({\rm Cat}))$.
\end{definition}

By duality, 
a double oplax monoidal functor
between duoidal categories
$(\mathcal{C},\boxtimes,1_{\boxtimes},\otimes,1_{\otimes})$
and
$(\mathcal{D},\boxtimes,1_{\boxtimes},\otimes,1_{\otimes})$
is a triple $(F,\varphi,\psi)$
such that
$(F^{\rm op},\psi^{\rm op},\varphi^{\rm op})$
is a double lax monoidal functor
from 
$(\mathcal{C}^{\rm op},\otimes,1_{\otimes},\boxtimes,1_{\boxtimes})$
to
$(\mathcal{D}^{\rm op},\otimes,1_{\otimes},\boxtimes,1_{\boxtimes})$
(see \cite[Definition~6.55]{Aguiar-Mahajan} for details).

\section{Monoidal $\infty$-categories}
\label{section:monoidal-infinity-categories}

In this section we recall the definition
of monoidal $\infty$-categories and strong monoidal functors
as monoid objects and morphisms between them
in the $\infty$-category of (small) $\infty$-categories.
Furthermore, 
we recall the definition of lax monoidal 
functors between monoidal $\infty$-categories
by regarding them as morphisms of 
nonsymmetric $\infty$-operads.
By reformulating monoidal $\infty$-categories 
in terms of Cartesian fibrations,
we introduce oplax monoidal functors.  
We think that results in this section
are well known %to experts
and do not claim any originality on them.

\subsection{Monoidal $\infty$-categories 
and strong monoidal functors}
\label{subsection:monoidal-infinity-categories-monoidal-functors}

In this subsection we recall the definitions 
of monoidal $\infty$-categories and strong monoidal functors
as monoid objects and morphisms between them
in the $\infty$-category
$\cat$ of (small) $\infty$-categories
(cf.~\cite[Definition~2.4.2.1]{Lurie2} and
\cite[Definition~3.5.1]{Gepner-Haugseng}).
%and
%\cite[Definition~4.14]{Groth}).
%(cf.~\cite[Proposition~3.5.4]{Gepner-Haugseng}).
%(cf.~\cite[Definition~2.2.13 and Proposition~3.5.4]{Gepner-Haugseng}).

Let $\Delta$ be the simplicial indexing category. 
A morphism $f: [n]\to [m]$ in $\Delta$ is said to be
inert if $f$ is an injection given by $[n]=\{0,1,\ldots,n\}\simeq
\{i,i+1,\ldots,i+n\}\hookrightarrow
\{0,1,\ldots,m\}=[m]$ for some $0\le i\le m-n$.
We denote by $\rho^{i}:[1]\to [n]$
the inert morphism in $\Delta$ given by
$[1]=\{0,1\}\simeq\{i-1,i\}\hookrightarrow
\{0,1,\ldots,n\}=[n]$ for $1\le i\le n$.

%A monoidal $\infty$-category
%is a monoid object in the $\infty$-category 
%$\cat$ of $\infty$-categories.
We recall the definition of monoid objects in 
an $\infty$-category with finite products
(cf.~\cite[Definition~2.4.2.1]{Lurie2} and
\cite[Definition~3.5.1]{Gepner-Haugseng}).

\begin{definition}\rm
Let $\mathcal{X}$ be an $\infty$-category 
which admits finite products,
and let $M: \Delta^{\rm op}\to \mathcal{X}$ be
a simplicial object in $\mathcal{X}$.
The inert morphism $\rho^i$ induces 
a morphism
$M(\rho^i): M([n])\to M([1])$
for $1\le i\le n$.
Taking a product of $M(\rho^i)$,
we obtain a morphism
\[ \prod_{1\le i\le n}M(\rho^i):
   M([n])\longrightarrow \overbrace{M([1])\times\cdots\times M([1])}^n.\]
We call this morphism a Segal morphism.
A monoid object in $\mathcal{X}$
is a simplicial object $M: \Delta^{\rm op}\to \mathcal{X}$
such that the Segal morphisms
are equivalences in $\mathcal{X}$ for all $[n]\in\Delta^{\rm op}$.
A morphism between monoid objects 
is a morphism of simplicial objects in $\mathcal{X}$.
We denote by
\[ \monno{\mathcal{X}} \]
the $\infty$-category of monoid objects in $\mathcal{X}$
and morphisms between them.
Note that $\monno{\mathcal{X}}$
is a full subcategory
of the $\infty$-category
${\rm Fun}(\Delta^{\rm op},\mathcal{X})$
of simplicial objects in $\mathcal{X}$.
\end{definition}

Let $\cat$ be the $\infty$-category
of (small) $\infty$-categories.
Since $\cat$ has finite products,
we can consider monoid objects in the $\infty$-category $\cat$.

\begin{definition}\rm
A monoidal $\infty$-category is a monoid
object in $\cat$,
and a strong monoidal functor 
between monoidal $\infty$-categories
is a morphism of monoid objects
in $\cat$.
The $\infty$-category
\[ \monno{\cat} \]
is the $\infty$-category
of monoidal $\infty$-categories and 
strong monoidal functors. 
\end{definition}

\subsection{Monoidal $\infty$-categories as coCartesian fibrations}
\label{subsection:monoidal-coCartesian-fibrations}

In this subsection we recall another formulation
of monoidal $\infty$-categories as coCartesian 
fibrations of nonsymmetric $\infty$-operads
(cf.~\cite[Example~2.4.2.4]{Lurie2}
and \cite[Definition~4.14]{Groth}).
%Using this formulation,
%we recall the definition of lax monoidal functors
%between monoidal $\infty$-categories. 

First, we recall the straightening and unstraightening
equivalence by Lurie \cite[Theorem~3.2.0.1]{Lurie1}.
Let $S$ be an $\infty$-category,
and let $\overcat{S}$
be the slice category of $\cat$
over $S$.
We denote by $\cocart{S^{\sharp}}$
the subcategory of $\overcat{S}$
consisting of coCartesian fibrations over $S$ and 
functors over $S$ which preserve coCartesian morphisms.
By Lurie's straightening and unstraightening
functors \cite[\S3.2]{Lurie1},
there is an equivalence 
\[ {\rm Fun}(S,\cat)\simeq \cocart{S^{\sharp}} \]
of $\infty$-categories.

Thus, monoid objects of $\cat$
corresponds to coCartesian fibrations
over $\Delta^{\rm op}$
with some properties 
under the equivalence between
${\rm Fun}(\Delta^{\rm op},\cat)$
and $\cocart{(\Delta^{\rm op})^{\sharp}}$.
We shall describe coCartesian fibrations 
over $\Delta^{\rm op}$ corresponding to
monoidal $\infty$-categories.

Let $p: X\to \Delta^{\rm op}$ be a coCartesian fibration.
We denote by $X_{[n]}$ the fiber of $p$ at $[n]\in\Delta^{\rm op}$.
The inert morphism $\rho^i$ of $\Delta$ induces a functor
$(\rho^i)_!: X_{[n]}\to X_{[1]}$ of $\infty$-categories.
Taking a product of $(\rho^i)_!$ for $1\le i\le n$,
we obtain a functor
\[ \prod_{1\le i\le n} (\rho^i)_!: X_{[n]}\longrightarrow
   \overbrace{X_{[1]}\times \cdots\times X_{[1]}}^n, \] 
which is also called a Segal morphism.
A coCartesian fibration $p: X\to\Delta^{\rm op}$
corresponds to a monoidal $\infty$-category
if and only if the Segal morphisms 
are equivalences for all $[n]\in\Delta^{\rm op}$.
We also call such a coCartesian fibration 
a monoidal $\infty$-category.

Next, we shall describe strong monoidal functors
between monoidal $\infty$-categories
in terms of coCartesian fibrations.
Let $p: X\to \Delta^{\rm op}$ and
$q: Y\to \Delta^{\rm op}$
be coCartesian fibrations
which are monoidal $\infty$-categories.
We consider a functor
$f: X\to Y$ 
of $\infty$-categories over $\Delta^{\rm op}$.
%\[ \xymatrix{
%    X \ar[rr]^{f} \ar[dr]_{p} & & Y \ar[dl]^{q}\\
%    & \Delta^{\rm op}. & 
%   }\]
The functor $f$ corresponds to a strong monoidal functor
if and only if $f$ carries $p$-coCartesian morphisms to
$q$-coCartesian morphisms.

We denote by
\[ \mathsf{Mon}(\cat) \]
the full subcategory of
$\cocart{(\Delta^{\rm op})^{\sharp}}$
spanned by monoidal $\infty$-categories.
By the above argument,
there is an equivalence of $\infty$-categories
\[ \monno{\cat}\simeq \mathsf{Mon}(\cat). \]

\begin{remark}\label{remark:monoidal-cat-cartesian-fibration}\rm
Let $\cart{\Delta^{\sharp}}$
be the subcategory of $\overcat{\Delta}$
consisting of Cartesian fibrations over $\Delta$
and functors over $\Delta$
which preserve Cartesian morphisms.
By Lurie's straightening and unstraightening
functors \cite[\S3.2]{Lurie1}, 
there is an equivalence
${\rm Fun}(\Delta^{\rm op},\cat)\simeq
   \cart{\Delta^{\sharp}}$
of $\infty$-categories.
Thus, we can also describe monoidal $\infty$-categories
in terms of Cartesian fibrations.
We say that a Cartesian fibration $p: X\to\Delta$
is a monoidal $\infty$-category 
if the Segal morphisms 
$\prod_{1\le i\le n}(\rho^i)^!:
X_{[n]}\to X_{[1]}\times\cdots\times X_{[1]}$
are equivalences
for all $[n]\in\Delta$.
We define 
\[ \mathsf{Mon}'(\cat) \]
to be the full subcategory of $\cart{\Delta^{\sharp}}$
spanned by monoidal $\infty$-categories.

Let $(-)^{\rm op}: \cat\to\cat$ be the functor
which assigns to an $\infty$-category $\mathcal{C}$
its opposite $\infty$-category $\mathcal{C}^{\rm op}$.
The functor $(-)^{\rm op}$ induces an equivalence
between $\cart{\Delta^{\sharp}}$
and $\cocart{(\Delta^{\rm op})^{\sharp}}$.
Furthermore,
this equivalence respects the Segal conditions.
Hence $(-)^{\rm op}$ induces an equivalence of $\infty$-categories
\[ (-)^{\rm op}: \mathsf{Mon}'(\cat)\simeq \mathsf{Mon}(\cat). \]
\end{remark}

\subsection{Lax and oplax monoidal functors}
\label{subsection:oplax-monoidal-functors}

In this subsection we recall the definition of
lax monoidal functors between monoidal $\infty$-categories
as morphisms of coCartesian fibrations of 
nonsymmetric $\infty$-operads
(cf.~\cite[Definition~2.1.2.7]{Lurie2}
and \cite[Definition~4.23]{Groth}).
By duality,
we introduce oplax monoidal functors
by using the description of monoidal $\infty$-categories
as Cartesian fibrations.

First, we recall the definition
of lax monoidal functors 
between monoidal $\infty$-categories.
Let $p: X\to \Delta^{\rm op}$
and $q: Y\to\Delta^{\rm op}$
be objects of $\mathsf{Mon}(\cat)$.
We consider a functor
$f: X\to Y$ of $\infty$-categories
over $\Delta^{\rm op}$.
%\[ \xymatrix{
%    X \ar[rr]^f\ar[dr]_p&& Y\ar[dl]^q\\
%    & \Delta^{\rm op}. & \\
%} \]

\begin{definition}
[{cf.~\cite[Definition~2.1.2.7]{Lurie2} and
 \cite[Definition~4.23]{Groth}}]
\rm
We say that $f$ is a lax monoidal functor
if $f$ carries $p$-coCartesian morphisms
over inert morphisms of $\Delta^{\rm op}$
to $q$-coCartesian morphisms.
%We denote by $\moncat{}$
%the $\infty$-category of monoidal $\infty$-categories
%and strong monoidal functors.
%We can regard $\moncat{}$
%as a full subcategory of $(\cocart{\Delta^{\rm op}})^{\rm strict}$.
We define an $\infty$-category 
\[ \laxmoncat \]
to be the subcategory of 
$\overcat{\Delta^{\rm op}}$
%$\cocart{(\Delta^{\rm op})^{\sharp}}$ 
consisting of monoidal $\infty$-categories and
lax monoidal functors.
%Note that we can regard $\laxmoncat{}$ as 
%a subcategory of $\cocart{\Delta^{\rm op}}$.
\end{definition}

For a monoidal $\infty$-category
$p: X\to\Delta^{\rm op}$,
we call the fiber $X_{[1]}$ at $[1]\in\Delta^{\rm op}$
the underlying $\infty$-category of 
the monoidal $\infty$-category.
Assigning to a monoidal $\infty$-category 
its underlying $\infty$-category,
we obtain a functor
\[ {\rm ev}_{[1]}: \laxmoncat{}\longrightarrow \cat.\]

\begin{remark}\rm
A monoidal $\infty$-category can be identified
with a non-symmetric $\infty$-operad 
$p: X\to \Delta^{\rm op}$
such that $p$ is a coCartesian fibration.
A lax monoidal functor between monoidal
$\infty$-categories is a map
of non-symmetric $\infty$-operads.
Thus, the $\infty$-category
$\laxmoncat$ is a full subcategory of $\nonoperad$
spanned by monoidal $\infty$-categories,
where $\nonoperad$ is the $\infty$-category 
of non-symmetric $\infty$-operads.
\end{remark}

Next,
we shall define an oplax monoidal functor
between monoidal $\infty$-categories.
For this purpose,
it is convenient for us to 
regard monoidal $\infty$-categories as 
Cartesian fibrations.
Let $p:X\to\Delta$ and $q: Y\to\Delta$
be objects of $\mathsf{Mon}'(\cat)$.
We consider a functor
$f: X\to Y$ of $\infty$-categories
over $\Delta$.
%\[ \xymatrix{
%    X \ar[rr]^f\ar[dr]_p&& Y\ar[dl]^q\\
%    & \Delta. & \\
%} \]

\begin{definition}\rm
We say that $f$ is an oplax monoidal functor
if $f$ carries $p$-Cartesian morphisms
over inert morphisms of $\Delta$
to $q$-Cartesian morphisms.
We define an $\infty$-category 
\[ \oplaxmoncat \]
to be the subcategory of 
$\overcat{\Delta}$
%$\cart{\Delta^{\sharp}}$ 
consisting of monoidal $\infty$-categories and
oplax monoidal functors.
\end{definition}
For an object $q: Y\to\Delta$ of
$\oplaxmoncat$,
we say that the fiber $Y_{[1]}$ at $[1]\in\Delta$
is the underlying
$\infty$-category of the monoidal $\infty$-category.
Assigning the underlying $\infty$-category
to a monoidal $\infty$-category, 
we obtain a functor
\[ {\rm ev}_{[1]}: \oplaxmoncat{}\longrightarrow \cat.\]

Finally,
we shall describe a relationship between $\oplaxmoncat{}$
and $\laxmoncat{}$
under the functor $(-)^{\rm op}$.
Recall that  
$(-)^{\rm op}: \cat\to \cat$
is the functor which assigns to an $\infty$-category 
$\mathcal{C}$
its opposite $\infty$-category $\mathcal{C}^{\rm op}$,
and that  
$(-)^{\rm op}$ induces an equivalence
between $\mathsf{Mon}'(\cat)$ and $\mathsf{Mon}(\cat)$.
%For a Cartesian fibration $p: X\to \Delta$,
%we obtain a coCartesian fibration
%$R(p): X^{\rm op}\to\Delta^{\rm op}$.
%by applying the functor $R:\cat\to\cat$,
%which assigns to an $\infty$-category 
%its opposite $\infty$-category.
%We observe that 
%$p: X\to \Delta$
%is an object of $\mathsf{Mon}'(\cat)$
%if only if $R(p): X^{\rm op}\to \Delta^{\rm op}$ 
%is an object of $\mathsf{Mon}(\cat)$.

Let $f: X\to Y$
be a morphism in $\overcat{\Delta}$
between monoidal $\infty$-categories. 
We observe that 
$f$ is an oplax monoidal functor
if and only if 
$f^{\rm op}$
%: X^{\rm op}\to Y{}^{\rm op}$ 
%in $\cart{\Delta}$
is a lax monoidal functor.
%between monoidal $\infty$-categoreis.
Hence we obtain a functor
\[ (-)^{\rm op}: \oplaxmoncat{}\longrightarrow \laxmoncat{}\]  
by assigning $p^{\rm op}: X^{\rm op}\to \Delta^{\rm op}$
to a monoidal $\infty$-category $p: X\to \Delta$.
We note that the functor
$(-)^{\rm op}:\oplaxmoncat{}\to\laxmoncat{}$ 
is an equivalence of $\infty$-categories,
%\begin{lemma}
and that there is a pullback diagram
\[ \begin{array}{ccc}
    \oplaxmoncat{} & \stackrel{(-)^{\rm op}}{\longrightarrow} &
    \laxmoncat{}\\[2mm]
    \mbox{$\scriptstyle {\rm ev}_{[1]}$}
    \bigg\downarrow\phantom{\mbox{$\scriptstyle {\rm ev}_{[1]}$}}
    & & 
    \phantom{\mbox{$\scriptstyle {\rm ev}_{[1]}$}}\bigg\downarrow
    \mbox{$\scriptstyle {\rm ev}_{[1]}$}\\[2mm]
    \cat & \stackrel{(-)^{\rm op}}{\longrightarrow} & \cat 
   \end{array}\]
in the $\infty$-category $\wcat$
of (large) $\infty$-categories.
%\end{lemma}

%\input{lax-monoidal-functors}
\section{Mixed fibrations}
\label{section:mixed-fibrations}

In this section we
introduce mixed fibrations
and study their properties.
We remark that
mixed fibrations are also studied
in \cite{HHLN1,HHLN2}
by name of curved orthofibration,
and in \cite{Stefanich}
by lax two-sided fibration. 
In \S\ref{subsection:Functor-infinity-categories}
we study monoid objects of the $\infty$-category
of a slice category $\overcat{T}$.
We consider (op)lax monoidal functors
between monoid objects of $\overcat{T}$.
In \S\ref{subsection:mixed-fibration},
motivated by 
the results in \S\ref{subsection:Functor-infinity-categories},
we introduce mixed fibrations.
Although mixed fibrations
have an obvious duality by definition,
we also give another asymmetric description 
of mixed fibrations.
%In \S\ref{subsection:marked-mixed-fibrations}
We consider mixed fibrations over marked simplicial sets
and call them marked mixed fibrations. 
By considering morphisms between marked mixed fibrations
which preserve (co)Cartesian morphisms over marked edges,
we introduce an $\infty$-category of marked mixed fibrations.

\subsection{Monoid objects of $\overcat{T}$
and (op)lax monoidal functors}
\label{subsection:Functor-infinity-categories}

%\subsection{$\infty$-categories of functors}

Let $T$ be an $\infty$-category.
The $\infty$-category $\overcat{T}$
of (small) $\infty$-categories over $T$ 
has finite products.
We consider a subcategory
$\mathcal{C}$ of $\overcat{T}$
which is closed under finite products
and equivalences
in $\overcat{T}$.
%{\color{red} Furthermore, we assume that
%the inclustion functor
%$\mathcal{C}\hookrightarrow\overcat{T}$
%is conservative.}
In this subsection we shall introduce
an $\infty$-category ${\rm Mon}^{\rm lax}(\mathcal{C})$ 
of monoid objects in $\mathcal{C}$ and lax monoidal functors.
By duality, we also introduce
an $\infty$-category ${\rm Mon}^{\rm oplax}(\mathcal{C})$
of monoid objects in $\mathcal{C}$ and oplax monoidal functors.

First, we consider
a description of the $\infty$-category
${\rm Fun}(S,\overcat{T})$ of functors
from an $\infty$-category $S$ to $\overcat{T}$
in terms of fibrations over $S\times T$.
Let $P: S\to \overcat{T}$ be a functor of $\infty$-categories.
The composite of $P$
with the projection 
$\overcat{T}\to\cat$
gives rise to a coCartesian fibration
$p_S: X\to S$
under the straightening and unstraightening
equivalence \cite[Theorem~3.2.0.1]{Lurie1}.
Furthermore,
we obtain a commutative diagram 
\[ \xymatrix{
    X \ar[rr]^{p}\ar[dr]_{p_S}    && S\times T\ar[dl]^{\pi_S} \\
    & S,   & \\
    }\]
where $\pi_S$ is the projection and
$p$ carries $p_S$-coCartesian morphisms to
$\pi_S$-coCartesian morphisms.
We may assume that
$p$ is a categorical fibration
by decomposing $p$
into a composite of a categorical equivalence
and a categorical fibration.

A morphism $F: P\to Q$ in ${\rm Fun}(S,\overcat{T})$
gives rise to a commutative diagram 
\[ \xymatrix{
    X \ar[rr]^{f}\ar[dr]_{p}    && Y\ar[dl]^{q} \\
    & S\times T,   & \\
    }\]
where $p$ and $q$ are categorical fibrations 
corresponding to $P$ and $Q$,
respectively, 
and $f$ sends $p_S$-coCartesian morphisms
to $q_S$-coCartesian morphisms.

Next, we consider an $\infty$-category
${\rm Fun}(S,\mathcal{C})$,
where $\mathcal{C}$ is a subcategory
of $\overcat{T}$
which is closed under 
%finite products and 
equivalences in $\overcat{T}$.
We shall give a similar description
of ${\rm Fun}(S,\mathcal{C})$ in terms
of fibrations over $S\times T$ as above.

For an object $P\in {\rm Fun}(S,\mathcal{C})$,
the composite of $P$ with 
the inclusion functor $\mathcal{C}\hookrightarrow\overcat{T}$
induces a categorical fibration $p: X\to S\times T$
such that $p_S=\pi_S\circ p$ is a coCartesian fibration
and $p$ preserves coCartesian morphisms as above.
We denote by $X_s$ the fiber of
$p_S: X\to S$ at $s\in S$.
By restriction of $p$,
we obtain a functor $p_s: X_s\to \{s\}\times T\simeq T$,
which is equivalent to $P(s)$ in $\mathcal{C}$
for any $s\in S$.
Furthermore,
a morphism $e:s\to s'$ in $S$
induces a functor $p_{e}: X_s\to X_{s'}$.
Since $p$ preserves coCartesian morphisms,
we obtain a commutative diagram
\[ \xymatrix{
    X_s \ar[rr]^{p_{e}}\ar[dr]_{p_s}    && 
    X_{s'}\ar[dl]^{p_{s'}} \\
    & T,   & \\
    }\]
which is equivalent to the morphism
$P(e): P(s)\to P(s')$ in $\mathcal{C}$.

Let $F: P\to Q$ be a morphism 
in ${\rm Fun}(S,\mathcal{C})$.
We have categorical fibrations $p: X\to S\times T$ and 
$q: Y\to S\times T$ which correspond to $P$ and $Q$,
respectively.
Regarding $F$ as a morphism 
in ${\rm Fun}(S,\overcat{T})$
by using the inclusion functor
$\mathcal{C}\hookrightarrow \overcat{T}$,
we obtain a functor 
$f: X\to Y$ over $S\times T$
such that 
$f$ carries $p_S$-coCartesian morphisms
to $q_S$-coCartesian morphisms.
By restriction of $f$,
we obtain the following commutative diagram 
\[ \xymatrix{
    X_s \ar[rr]^{f_s}\ar[dr]_{p_s}    && Y_s\ar[dl]^{q_s} \\
    & T,   & \\
    }\]
which is equivalent
to $F(s): P(s)\to Q(s)$ in $\mathcal{C}$
for any $s\in S$.

Motivated by the description of
${\rm Fun}(S,\mathcal{C})$ in terms of
fibrations over $S\times T$ as above,
we introduce an $\infty$-category
$\mathsf{Fun}(S,\mathcal{C})$ which is equivalent 
to ${\rm Fun}(S,\mathcal{C})$.

\begin{definition}\rm\label{def:mathsf-fun-S-C}
Let $S$ and $T$ be $\infty$-categories, and
let $\mathcal{C}$ be a subcategory of
$\overcat{T}$ which is closed under %finite products and 
equivalences in $\overcat{T}$.      
We define a subcategory $\mathsf{Fun}(S,\mathcal{C})$  
of $\overcat{S\times T}$ as follows.
The objects of $\mathsf{Fun}(S,\mathcal{C})$
are categorical fibrations $p: X\to S\times T$ such that
$p_S= \pi_S\circ p:X\to S$ is a coCartesian fibration,
$p$ carries $p_S$-coCartesian morphisms
to $\pi_S$-coCartesian morphisms,
the restriction $p_s: X_s\to T$ is an object of $\mathcal{C}$
for any $s\in S$,
and the induced functor
$p_e: X_s\to X_{s'}$ over $T$ 
is a morphism in $\mathcal{C}$
for any morphism $e: s\to s'$ in $S$.
The morphisms of $\mathsf{Fun}(S,\mathcal{C})$ between
$p: X\to S\times T$ and $q:Y\to S\times T$
are functors $f: X\to Y$ over $S\times T$
such that
$f$ carries $p_S$-coCartesian morphisms 
to $q_S$-coCartesian morphisms
and 
the induced functor $f_s: X_s\to Y_s$ over $T$
is a morphism in $\mathcal{C}$ for any $s\in S$.
\end{definition}

%By the above argument,
%we obtain the following proposition.

\begin{proposition}\label{prop:description-fun-S-general-C}
There is an equivalence of $\infty$-categories
\[ \mathsf{Fun}(S,\mathcal{C})\simeq {\rm Fun}(S,\mathcal{C}). \]
\end{proposition}

\begin{remark}\rm
Let $S$ and $T$ be $\infty$-categories, and
let $\mathcal{C}$ be a subcategory of
$\overcat{T}$ which is closed under 
%finite products and 
equivalences in $\overcat{T}$.      
By duality,
we define a subcategory $\mathsf{Fun}'(S^{\rm op},\mathcal{C})$  
of $\overcat{S\times T}$ as follows.
The objects of $\mathsf{Fun}'(S^{\rm op},\mathcal{C})$
are categorical fibrations $p: X\to S\times T$ such that
$p_S= \pi_S\circ p:X\to S$ is a Cartesian fibration,
$p$ carries $p_S$-Cartesian morphisms
to $\pi_S$-Cartesian morphisms,
the restriction $p_s: X_s\to T$ is an object of $\mathcal{C}$
for any $s\in S$,
and the induced functor
$p_e: X_s\to X_{s'}$ over $T$ 
is a morphism in $\mathcal{C}$
for any morphism $e: s'\to s$ in $S$.
The morphisms of $\mathsf{Fun}'(S^{\rm op},\mathcal{C})$ between
$p: X\to S\times T$ and $q:Y\to S\times T$
are functors $f: X\to Y$ over $S\times T$
such that
$f$ carries $p_S$-Cartesian morphisms 
to $q_S$-Cartesian morphisms
and 
the induced functor $f_s: X_s\to Y_s$ over $T$
is a morphism in $\mathcal{C}$ for any $s\in S$.
By the same argument as above,
we see that 
there is an equivalence of $\infty$-categories
\[ \mathsf{Fun}'(S^{\rm op},\mathcal{C})\simeq {
   \rm Fun}(S^{\rm op},\mathcal{C}). \]
\end{remark}

Next,
we shall consider monoid objects in 
a subcategory $\mathcal{C}$ of $\overcat{T}$. 
We assume that $\mathcal{C}$ is 
closed under finite products and equivalences.
Since $\mathcal{C}$ is closed under finite products,
we can consider the $\infty$-category
${\rm Mon}(\mathcal{C})$ of monoid objects
in $\mathcal{C}$.
The $\infty$-category
${\rm Mon}(\mathcal{C})$
is a full subcategory of 
${\rm Fun}(\Delta^{\rm op},\mathcal{C})$
and we have an equivalence between
${\rm Fun}(\Delta^{\rm op},\mathcal{C})$
and $\mathsf{Fun}(\Delta^{\rm op},\mathcal{C})$
by Proposition~\ref{prop:description-fun-S-general-C}.
This motivates us to introduce an
$\infty$-category $\mathsf{Mon}(\mathcal{C})$ 
as follows.

\begin{definition}\label{definition:mathsfmon-C}\rm
We define an $\infty$-category
\[ \mathsf{Mon}(\mathcal{C}) \]
to be a full subcategory of 
$\mathsf{Fun}(\Delta^{\rm op},\mathcal{C})$ 
spanned by those objects $p: X\to \Delta^{\rm op}\times T$ 
of $\mathsf{Fun}(\Delta^{\rm op},\mathcal{C})$ 
such that 
the induced morphisms
\[ X_{[n]}\longrightarrow \overbrace{X_{[1]}\times_T\cdots
          \times_TX_{[1]}}^n \]
are equivalences in $\mathcal{C}$
for all $[n]\in\Delta^{\rm op}$. 
\end{definition}

\begin{remark}\label{remark:mathsfMon-primes-C}\rm
By duality, we define an $\infty$-category $\mathsf{Mon}'(\mathcal{C})$ 
to be a full subcategory of 
$\mathsf{Fun}'(\Delta,\mathcal{C})$ 
spanned by those objects $p: X\to \Delta\times T$ 
of $\mathsf{Fun}'(\Delta,\mathcal{C})$ 
such that the induced morphisms
$X_{[n]}\to X_{[1]}\times_T\cdots\times_TX_{[1]}$
are equivalences in $\mathcal{C}$
for all $[n]\in\Delta^{\rm op}$. 
\end{remark}

Using Proposition~\ref{prop:description-fun-S-general-C},
we can verify that
a monoid object of $\mathcal{C}$ corresponds
to an object of $\mathsf{Mon}(\mathcal{C})$
and a morphism between monoid objects in $\mathcal{C}$ 
corresponds to a morphism of $\mathsf{Mon}(\mathcal{C})$. 
Thus, we obtain the following proposition.

\begin{proposition}\label{prop:moncal-monsf}
There is an equivalence of $\infty$-categories
\[ \mathsf{Mon}(\mathcal{C})\simeq \monno{\mathcal{C}}. \]
\end{proposition}

\begin{remark}\rm
By duality, there is an equivalence
$\mathsf{Mon}'(\mathcal{C})\simeq \monno{\mathcal{C}}$
of $\infty$-categories.
\end{remark}

\begin{remark}\rm
Let $p: X\to \Delta^{\rm op}\times T$
be an object of $\mathsf{Mon}(\mathcal{C})$.
%We denote by $X_t$
%the fiber of $\pi_T\circ p: X\to T$ at $t\in T$,
%and by $p_t: X_t\to\Delta^{\rm op}\times\{t\}\simeq\Delta^{\rm op}$
%the restriction of $p$.
The restriction $p_t: X_t\to \Delta^{\rm op}$ is a monoidal
$\infty$-category for each $t\in T$.
%Let $p: X\to \Delta^{\rm op}\times T$
%and $q: Y\to \Delta^{\rm op}\times T$
%be objects of $\mathsf{Mon}(\mathcal{C})$,
Let $f: X\to Y$ be a morphism of $\mathsf{Mon}(\mathcal{C})$.
Then $f$ induces a strong monoidal functor
$f_t: X_t\to Y_t$ for each $t\in T$.
\end{remark}

We shall consider a lax monoidal functor
between monoid objects of $\mathcal{C}$.

\begin{definition}\rm
\label{def:lax-mon-cat-C}
Let $f: X\to Y$ be a functor over
$\Delta^{\rm op}\times T$,
where $p: X\to \Delta^{\rm op}\times T$
and $q: Y\to \Delta^{\rm op}\times T$
are objects of $\mathsf{Mon}(\mathcal{C})$.
We say that $f$ is a lax monoidal functor 
if $f_{[n]}: X_{[n]}\to Y_{[n]}$
is a morphism of $\mathcal{C}$ for each $[n]\in\Delta^{\rm op}$,
and if $f$ carries $p_{\Delta^{\rm op}}$-coCartesian
morphisms over inert morphisms of $\Delta^{\rm op}$
to $q_{\Delta^{\rm op}}$-coCartesian morphisms.
We denote by
\[ \laxmon(\mathcal{C}) \]
the subcategory of $\cat{}_{/\Delta^{\rm op}\times T}$
consisting of objects of $\mathsf{Mon}(\mathcal{C})$
and lax monoidal functors between them.
\end{definition}

\begin{remark}\rm
%If $p: X\to \Delta^{\rm op}\times T$
%is an object of $\laxmon(\mathcal{C})$,
%then $p_t: X_t\to \Delta^{\rm op}$ are monoidal
%$\infty$-categories for all $t\in T$.
%Let $p: X\to \Delta^{\rm op}\times T$
%and $q: Y\to \Delta^{\rm op}\times T$
%be objects of $\laxmon{\mathcal{C}}$,
Let $f: X\to Y$ be a functor
of $\laxmon(\mathcal{C})$.
Then $f$ induces a lax monoidal functor
$f_t: X_t\to Y_t$ for each $t\in T$.
\end{remark}

By duality,
we shall define an $\infty$-category
$\oplaxmon(\mathcal{C})$
of monoid objects in $\mathcal{C}$
and oplax monoidal functors.

\begin{definition}\rm\label{def:description-oplaxmoncat}
Let $f: X\to Y$ be a functor over $\Delta\times T$,
where $p: X\to \Delta\times T$
and $q: Y\to \Delta\times T$ are objects of $\mathsf{Mon}'(\mathcal{C})$.
We say that $f$ is an oplax monoidal functor
if $f_{[n]}: X_{[n]}\to Y_{[n]}$
is a morphism of $\mathcal{C}$ for each $[n]\in\Delta$,
and if $f$ carries $p_{\Delta}$-Cartesian
morphisms over inert morphisms of $\Delta$
to $q_{\Delta}$-Cartesian morphisms.
We denote by 
\[ \oplaxmon(\mathcal{C}) \]
the subcategory of $\cat{}_{/\Delta\times T}$
consisting of objects of $\mathsf{Mon}'(\mathcal{C})$
and oplax monoidal functors between them.
\end{definition}

\begin{remark}\rm
Let $f: X\to Y$ be a functor
of $\oplaxmon(\mathcal{C})$.
Then $f$ induces an oplax monoidal functor
$f_t: X_t\to Y_t$ for each $t\in T$.
\end{remark}

\subsection{Mixed fibrations}
\label{subsection:mixed-fibration}

In this subsection, 
motivated by \S\ref{subsection:Functor-infinity-categories},
we introduce mixed fibrations and study their properties.
We remark that
mixed fibrations are also studied
in \cite{HHLN1,HHLN2}
by name of curved orthofibration,
and in \cite{Stefanich}
by lax two-sided fibration. 

Let 
$p: X\to S\times T$
be a categorical fibration,
where $S$ and $T$ are $\infty$-categories.
%Recall that $\pi_S: S\times T\to S$ 
%and $\pi_T:S\times T\to T$ are the projections.
%We consider the following conditions on $p$:
%
%Let $p: X\to S\times T$ be a categorical fibration.
We set $p_S=\pi_S\circ p$ and $p_T=\pi_T\circ p$,
where $\pi_S: S\times T\to S$ and 
$\pi_T: S\times T\to T$ are the projections.
We consider the following conditions on $p$:

\begin{condition}\rm\label{item:coCart-preserve-coCart}
The map $p_S$ is a coCartesian fibration and
the functor $p$ carries $p_S$-coCartesian morphisms
to $\pi_S$-coCartesian morphisms.
\end{condition}

\begin{condition}\rm\label{item:Cart-preserve-Cart}
The map $p_T$ is a Cartesian fibration and
the functor $p$ carries $p_T$-Cartesian morphisms
to $\pi_T$-Cartesian morphisms.
\end{condition}

\begin{definition}[{cf.~\cite[Definition~2.3.1]{HHLN1} and
\cite[Definition~2.1.1]{Stefanich}}]
\rm
Let $S$ and $T$ be $\infty$-categories.
We say that a categorical fibration 
$p: X\to S\times T$ is a mixed fibration over $(S,T)$
if it satisfies 
%Conditions~\ref{condition:mixed-coCartesian}
%and \ref{condition:mixed-Cartesian}.
Conditions~\ref{item:coCart-preserve-coCart}
and \ref{item:Cart-preserve-Cart}.
We define an $\infty$-category
\[ \mixed{(S,T)} \]
to be the full subcategory of $\overcat{S\times T}$
spanned by mixed fibrations over $(S,T)$,
and call it the $\infty$-category of mixed fibrations
over $(S,T)$.
\end{definition}

\begin{remark}\rm\label{remark:mixed-fib-R-duality}
%Recall that $(-)^{\rm op}: \cat\to\cat$ is the functor
%which assigns to an $\infty$-category
%its opposite $\infty$-category.
The functor $(-)^{\rm op}$ induces an equivalence
of $\infty$-categories between
$\overcat{S\times T}$ and $\overcat{T^{\rm op}\times S^{\rm op}}$.
By symmetry of the definition of mixed fibrations,
$(-)^{\rm op}$ induces an equivalence of $\infty$-categories 
\[ \mixed{(S,T)}\simeq \mixed{(T^{\rm op},S^{\rm op})}.  \]
\end{remark}

We consider characterizations of mixed fibrations.
We denote by $X_s$ the fiber of
the map $p_S: X\to S$ at $s\in S$,
and by $X_t$ the fiber of the map
$p_T: X\to T$ at $t\in T$. 
By restriction,
there are  maps
$p_s: X_s\to \{s\}\times T\simeq T$
and 
$p_t: X_t\to S\times\{t\}\simeq S$. 

%\begin{remark}\label{remark:p-cocart-in-Xt}\rm
%Let $p: X\to S\times T$ be a categorical fibration
%which satisfies Condition~\ref{item:coCart-preserve-coCart}.
%For any $x\in X_t\ (t\in T)$ and any morphism 
%$\overline{e}: p_S(x)\to s$ of $S$,
%we can take a $p_S$-coCartesian morphism $x\to x'$ 
%in $X_t$ over $\overline{e}$. 
%\end{remark}
%
%\begin{remark}\rm\label{remark:p-cart-in-Xs}
%Let $p: X\to S\times T$ be a categorical fibration
%which satisfies Condition~\ref{item:Cart-preserve-Cart}.
%By duality,
%for any $x\in X_s\ (s\in S)$ and 
%any morphism $\overline{e}: t\to p_T(x)$ of $T$,
%we can take a $p_T$-Cartesian morphism $x'\to x$ in 
%$X_s$ over $\overline{e}$.
%\end{remark}

%%\if0
\begin{lemma}\label{lemma:pt-cocart-eq-p-cocart}
%%Let $p: X\to S\times T$ be a mixed fibration
%%over $(S,T)$,
Let $p: X\to S\times T$ be a categorical fibration
and let $e: x\to x'$ be a morphism in $X_t$
for $t\in T$.
If $e$ is a $p_S$-coCartesian morphism,
then it is a $p_t$-coCartesian morphism.
Conversely,
if $e$ is a $p_t$-coCartesian morphism
and if $p: X\to S\times T$ satisfies 
Condition~\ref{item:coCart-preserve-coCart},
then $e$ is a $p_S$-coCartesian morphism. 
\end{lemma}

\proof
If $e$ is $p_S$-coCartesian,
then $e$ is $p$-coCartesian
by \cite[Proposition~2.2.1]{HHLN1}.
Since $p_t: X_t\to S$ is a pullback
of $p: X\to S\times T$ along 
$S\simeq S\times\{t\}\to S\times T$,
we see that $e$ is $p_t$-coCartesian.

Now,
assume that $e$ is $p_t$-coCartesian and
that $p$ satisfies Condition~\ref{item:coCart-preserve-coCart}.
By \cite[Corollary~2.2.2]{HHLN1},
there is a $p_S$-coCartesian morphism $e': x\to x''$ in $X_t$
over $p_S(e)$.
By the first part of the proof,
$e'$ is $p_t$-coCartesian.
This implies that 
$e$ and $e'$ are equivalent in $X_t$,
and hence $e$ is $p_S$-coCartesian.
\qed
%%\fi

\begin{remark}\rm\label{remark:p-cart-equivalent-pt-cart}
%Let $p: X\to S\times T$ be a mixed fibration
%over $(S,T)$, and 
Let $p: X\to S\times T$ be a categorical fibration
and let $e: x\to x'$ be a morphism in $X_s$
for $s\in S$.
By duality,
if $e$ is a $p_T$-Cartesian morphism,
then $e$ is a $p_s$-Cartesian morphism.
Conversely,
if $e$ is a $p_s$-Cartesian morphism
and if $p$ satisfies Condition~\ref{item:Cart-preserve-Cart},
then $e$ is a $p_T$-Cartesian morphism. 
\end{remark}
%\fi

We consider the following conditions:

\begin{condition}\rm\label{item:pt-coCart} 
The maps $p_t: X_t\to S$ are coCartesian fibrations
for all $t\in T$. 
\end{condition}

\begin{condition}\rm\label{item:qs-Cart}
The maps $p_s: X_s\to T$ are Cartesian fibrations
for all $s\in S$.
\end{condition}

We have the following characterizations
of mixed fibrations
by \cite[Proposition~2.3.3]{HHLN1}.

\begin{proposition}[{\cite[Proposition~2.3.3]{HHLN1}}]
\label{prop:characterization-mixed-fib}
The following conditions are
equivalent for $p: X\to S\times T$\mbox{\rm :}

\begin{enumerate}

\item[\mbox{\rm (1)}]
The map $p$ is a mixed fibration.

\item[\mbox{\rm (2)}]
The map $p$ satisfies
Conditions~\ref{item:coCart-preserve-coCart} and
\ref{item:qs-Cart}.

\item[\mbox{\rm (3)}]
The map $p$ satisfies
Conditions~\ref{item:Cart-preserve-Cart} and
\ref{item:pt-coCart}.
  
\end{enumerate}
\end{proposition}

%\input{marked-mixed-fibration}
%\subsection{Marked mixed fibrations}
%\label{subsection:marked-mixed-fibrations}

Now, 
we study properties of maps
between mixed fibrations.
For this purpose,
we introduce a notion of marked mixed fibrations.

First, we recall the definition of marked simplicial sets
(cf.~\cite[Definition~3.1.0.1]{Lurie1}).
A marked simplicial set is a pair $(S,E)$,
where $S$ is a simplicial set and
$E$ is a set of edges of $S$
which contains all degenerate edges.
%We denote by $S_E$ the pair $(S,E)$.
When $E$ is the set $S_1$ of all edges of $S$,
we write $S^{\sharp}=(S,S_1)$.
When $E$ is the set $s_0(S_0)$ of all degenerate edges of $S$,
we write $S^{\flat}=(S,s_0(S_0))$.

Let $(S,E)$ be a marked simplicial set,
where $S$ is an $\infty$-category.
We define an $\infty$-category 
\[ \cocart{(S,E)} \]
to be a subcategory of $\overcat{S}$
as follows.
The objects of $\cocart{(S,E)}$ 
are coCartesian fibrations over $S$.
The morphisms between coCartesian fibrations 
$p: X\to S$ and $q: Y\to S$
are functors $f: X\to Y$ over $S$
such that $f$ preserves coCartesian morphisms 
over $E$.

%{\color{red} Lurie's categorical pattern との関係について言及せよ。}
\begin{remark}\rm
Let $\mathfrak{P}=(E,T,\emptyset)$
be a categorical pattern on $S$
in the sense of \cite[Definition~B.0.19]{Lurie2},
where $T$ is the set of all $2$-simplices of $S$.
By \cite[Theorem~B.0.20]{Lurie2},
there exists a left proper combinatorial 
simplicial model structure on $({\rm Set}_{\Delta}^{+})_{/\mathfrak{P}}$.
We can regard the $\infty$-category $\cocart{(S,E)}$
as a full subcategory of the underlying $\infty$-category
of $({\rm Set}_{\Delta}^{+})_{/\mathfrak{P}}$.
\end{remark}

Similarly,
we define an $\infty$-category
$\cart{(S,E)}$ to be a subcategory
of $\overcat{S}$ consisting of
Cartesian fibrations over $S$
and functors over $S$ which
preserve Cartesian morphisms over $E$.

Let $S$ and $T$ be $\infty$-categories,
and let $p: X\to S\times T$
be a mixed fibration over $(S,T)$.
%functor of $\infty$-categories
%satisfying Conditions~\ref{item:coCart-preserve-coCart} 
%and \ref{item:Cart-preserve-Cart}.
For any morphism $e: s\to s'$ of $S$,
we have a commutative diagram
\[ \xymatrix{
    X_s \ar[rr]^{e_!} \ar[dr]_{p_s}& & X_{s'}\ar[dl]^{p_{s'}}\\
    & T. &\\
    }\]
Similarly, for any morphism $f: t'\to t$ of $T$,
we have a commutative diagram
\[ \xymatrix{
    X_{t} \ar[rr]^{f^!} \ar[dr]_{p_{t}}& & X_{t'}\ar[dl]^{p_{t'}}\\
    & S. &\\
    }\]
Suppose that we have sets $E$ and $F$ of
morphisms of $S$ and $T$, respectively,
which contain all degenerate ones.
We consider the following conditions on $p$:

\begin{condition}\rm\label{item:qs-cocart-inert-exchange}
For any $e: s\to s'$ in $E$,
the functor $e_!$ carries $p_s$-Cartesian morphisms
over $F$ to $p_{s'}$-Cartesian morphisms. 
\end{condition}

\begin{condition}\rm\label{item:pt-coCart-inert-exchange}
For any $f: t'\to t$ in $F$,
the functor $f^!$ carries $p_{t}$-coCartesian morphisms over $E$ 
to $p_{t'}$-coCartesian morphisms.
\end{condition}

\begin{proposition}\label{prop:cocart-cart-interchange}
Let $p: X\to S\times T$ be a mixed fibration
over $(S,T)$,
%functor of $\infty$-categories 
%satisfying Conditions~\ref{item:coCart-preserve-coCart} 
%and \ref{item:Cart-preserve-Cart}.
and let $E$ and $F$ be subsets of morphisms of $S$ and $T$,
respectively, which contain all degenerate ones.
The mixed fibration $p$ satisfies 
Condition~\ref{item:qs-cocart-inert-exchange} 
if and only if
$p$ satisfies Condition~\ref{item:pt-coCart-inert-exchange}.
\end{proposition}

\proof
We denote by $X_{(s,t)}$
the fiber of $p: X\to S\times T$ at $(s,t)\in S\times T$.
The proposition follows from the fact that
each of Conditions~\ref{item:qs-cocart-inert-exchange} 
and \ref{item:pt-coCart-inert-exchange} 
is equivalent to the condition that
there is a commutative diagram
\[ \begin{array}{ccc}
     X_{(s,t)} & \stackrel{e_!}{\longrightarrow} &
     X_{(s',t)}\\[1mm]
     \mbox{$\scriptstyle f^!$}\bigg\downarrow
     \phantom{\mbox{$\scriptstyle f^!$}} 
     & & 
     \phantom{\mbox{$\scriptstyle f^!$}}
     \bigg\downarrow\mbox{$\scriptstyle f^!$}\\[3mm]
     X_{(s,t')} & \stackrel{e_!}{\longrightarrow} &
     X_{(s',t')}\\
   \end{array}\]
for any $(e,f)\in E\times F$.
\qed

\bigskip

\begin{definition}\rm
We define a marked mixed fibration over $((S,E),(T,F))$
to be a mixed fibration $p: X\to S\times T$
over $(S,T)$ which satisfies
the equivalent Conditions 
\ref{item:qs-cocart-inert-exchange} and 
\ref{item:pt-coCart-inert-exchange}.
Let $p: X\to S\times T$ and $q: Y\to S\times T$
be marked mixed fibrations over $((S,E),(T,F))$.
A morphism $f$ of marked mixed fibrations over $((S,E),(T,F))$
between $p$ and $q$ is 
a morphism of mixed fibrations over $(S,T)$
such that $f$ carries $p_S$-coCartesian morphisms
over $E$ to $q_S$-coCartesian morphisms and
that $p_T$-Cartesian morphisms over $F$
to $q_T$-Cartesian morphisms. 
We define an $\infty$-category
\[ \mixed{((S,E),(T,F))} \]
to be the $\infty$-category 
of marked mixed fibrations over $((S,E),(T,F))$
and morphisms between them.
\end{definition}

\begin{remark}\rm
The notion of marked mixed fibration over $(S^{\sharp},T^{\sharp})$
is equivalent to that of orthofibration over
$T\times S$ in \cite{HHLN1} and of
two-sided fibration over $S\times T$
in \cite{Stefanich}.
\end{remark}
  
Next,
we consider a description of the $\infty$-category
${\rm Fun}(S,\cart{(T,F)})$
in terms of marked mixed fibrations.

\begin{proposition}\label{prop:description-fun-S-cart-T-pair-version}
There is an equivalence of $\infty$-categories   
\[ {\rm Fun}(S,\cart{(T,F)})\simeq
   \mixed{(S^{\sharp},(T,F))}. \]
\end{proposition}

\proof
By Proposition~\ref{prop:description-fun-S-general-C},
${\rm Fun}(S,\cart{(T,F)})$ is equivalent 
to $\mathsf{Fun}(S,\cart{(T,F)})$.
By Definition~\ref{def:mathsf-fun-S-C},
an object of $\mathsf{Fun}(S,\cart{(T,F)})$
is a categorical fibration $p: X\to S\times T$ 
such that $p$ satisfies 
Conditions~\ref{item:coCart-preserve-coCart},
$p_s: X_s\to T$ is a Cartesian fibration for any $s\in S$,
and 
$p_e: X_s\to X_{s'}$ over $T$ carries
$p_s$-Cartesian morphisms over $F$ to 
$p_{s'}$-Cartesian morphisms
for any morphism $e: s\to s'$ in $S$. 
%By Lemma~\ref{prop:S-cart-to-cart-fib},
By Proposition~\ref{prop:characterization-mixed-fib},
we see that 
$p$ is a mixed fibration over $(S,T)$.
Thus,
the objects of $\mathsf{Fun}(S,\cart{(T,F)})$
are the marked mixed fibrations $p: X\to S\times T$ 
over $(S^{\sharp},(T,F))$.
%\ref{item:Cart-preserve-Cart}, 
%which satisfies Condition~\ref{item:qs-cocart-inert-exchange}
%for $E=S_1$.

Let $p: X\to S\times T$ and $q: Y\to S\times T$
be objects of $\mathsf{Fun}(S,\cart{(T,F)})$.
A morphism between 
$p$ and $q$ in $\mathsf{Fun}(S,\cart{(T,F)})$
is a functor $f: X\to Y$ over $S\times T$ such that
$f$ carries $p_S$-coCartesian morphisms 
to $q_S$-coCartesian morphisms and 
the functor $f_s: X_s\to Y_s$ over $T$
carries $p_s$-Cartesian morphisms over $F$
to $q_s$-Cartesian morphisms for any $s\in S$.
By Remark~\ref{remark:p-cart-equivalent-pt-cart},
we see that
$f$ carries  $p_T$-Cartesian morphisms over $F$
to $q_T$-Cartesian morphisms.
%We shall show that $f$ carries  $p_T$-Cartesian morphisms over $F$
%to $q_T$-Cartesian morphisms.
%Let $\alpha: x\to x'$ be a $p_T$-Cartesian morphism
%over $F$.
%Since $p$ is a mixed fibration,
%we may assume taht $\alpha$ lies in $X_s$
%for some $s\in S$.
%Then, $\alpha$ is a $p_s$-Cartesian morphism over $F$.
%Hence, $f(\alpha)$ is a $q_s$-Cartesian morphism.
Thus, 
the morphisms between 
$p$ and $q$ in $\mathsf{Fun}(S,\cart{(T,F)})$
are the morphisms of mixed fibrations
over $(S^{\sharp},(T,F))$.
%By Proposition~\ref{prop:description-fun-S-general-C},
%the $\infty$-category 
%${\rm Fun}(S,(\cart{T})_E)$ is equivalent
%to $\mathsf{Fun}(S,(\cart{T})_E)$.
%
%Now, we consider another description of 
%objects of ${\rm Fun}(S,(\cart{T})^E)$.
%Let $f: X\to S\times T$ be an object
%of $\mathfrak{Fun}(S,\cart{T})$.
%This is equivalent to assume that 
%$f$ satisfies Conditions~\ref{item:coCart-preserve-coCart} 
%and \ref{item:Cart-preserve-Cart}.
%For any $\alpha: t\to t'$ of $T$,
%we have a commutative diagram
%\[ \xymatrix{
%     X_{t'}\ar[rr]^{\alpha^!}\ar[dr]_{p_{t'}}&&X_{t}\ar[dl]^{p_t}\\
%     & S.  &\\
%     }\]
%Hence we obtain the following proposition.
\qed

\begin{remark}\rm
By duality,
there is an equivalence of $\infty$-categories
\[ {\rm Fun}(T^{\rm op},\cocart{(S,E)})\simeq
   \mixed{((S,E),T^{\sharp})}.\]
\end{remark}

\section{Duoidal $\infty$-categories}
\label{section:duoidal-infty-categories}

%{\color{red}
In \S\ref{section:review-duoidal}
we recalled the notion of
duoidal category 
in the classical setting.
%}
In this section we introduce a duoidal category
in the setting of $\infty$-categories
and call it a duoidal $\infty$-category.
%{\color{red}
In \S\ref{section:review-duoidal}
we also recalled that
there are three kinds of functors between duoidal categories.
%}
We introduce corresponding functors between
duoidal $\infty$-categories.
According to which kind of functors of duoidal $\infty$-categories
we take,
we will give three formulations
of $\infty$-category of duoidal $\infty$-categories.

%\newpage

\subsection{Duoidal $\infty$-categories}

In this subsection we introduce duoidal $\infty$-categories
which are analogues of duoidal categories 
in the setting of $\infty$-categories.

We recall that $\oplaxmoncat{}$
is the $\infty$-category of
monoidal $\infty$-categories and oplax monoidal functors.
%which is a subcategory of $\cart{\Delta}$.
The objects of $\oplaxmoncat{}$
are Cartesian fibrations
$p: X\to \Delta$ such that
$p^{\rm op}: X^{\rm op}\to \Delta^{\rm op}$
is a monoidal $\infty$-category.
The morphisms between 
objects $X\to\Delta$ and $Y\to\Delta$
are functors $f:X\to Y$ over $\Delta$
such that $f$ preserves Cartesian morphisms
over inert morphisms of $\Delta$.
Let $\Delta_{\rm int}$ be the set of
all inert morphisms of $\Delta$.
We set $\Delta^{\natural}=(\Delta,\Delta_{\rm int})$.
Then $\oplaxmoncat{}$ is a full subcategory
of $\cart{\Delta^{\natural}}$
spanned by monoidal $\infty$-categories.

We can verify that $\oplaxmoncat{}$ has finite products
as follows.
For objects 
$p: X\to \Delta$ and $q: Y\to \Delta$
of $\oplaxmoncat{}$,
the fiber product 
$p\times_{\Delta}q:
X\times_{\Delta}Y\to\Delta$ 
is a product of $p$ and $q$,
and the identity functor
$\Delta\to\Delta$
is a final object in $\oplaxmoncat{}$. 
Thus, we can consider
monoid objects in $\oplaxmoncat{}$.

\begin{definition}
\label{definition:duoidal-infty-category-monoid-oplax}
\rm
A duoidal $\infty$-category
is a monoid object in $\oplaxmoncat{}$.
\end{definition}

\begin{example}\rm
Braided monoidal categories
can be regarded as duoidal categories.
We consider a generalization of this fact in $\infty$-categories.
Let $\mathbb{E}_2^{\otimes}$ be the little $2$-cubes operad.
We will show that an $\mathbb{E}_2$-monoidal
$\infty$-category determines a duoidal $\infty$-category.
Suppose that 
$m:\mathbb{E}_2^{\otimes}\to\cat$ 
is an $\mathbb{E}_2$-monoid object of $\cat$.
Using the functor ${\rm Cut}: \Delta^{\rm op}\to {\rm Assoc}^{\otimes}$
of \cite[Construction~4.1.2.9]{Lurie2},
we obtain a functor
$d:\Delta^{\rm op}\times \Delta^{\rm op}
{\to}{\rm Assoc}^{\otimes}\times{\rm Assoc}^{\otimes}
\simeq 
\mathbb{E}_1^{\otimes}\times\mathbb{E}_1^{\otimes}
\to
\mathbb{E}_2^{\otimes}$.
%\hookrightarrow \mathbb{E}_k^{\otimes}$.
By composing $m$ with $d$,
we obtain a functor
$d\circ m: \Delta^{\rm op}\times\Delta^{\rm op}\to\cat$,
which determines a functor
$\Delta^{\rm op}\to {\rm Mon}(\cat)\subset
{\rm Fun}(\Delta^{\rm op},\cat)$. 
By regarding it as a functor
$\Delta^{\rm op}\to {\rm Mon}(\cat)\simeq
\mathsf{Mon}'(\cat)\subset\mathsf{Mon}^{\rm oplax}(\cat)$,
we obtain a duoidal $\infty$-category.
\end{example}

%\bigskip

We give another description
of duoidal $\infty$-categories
in terms of mixed fibrations.
Let $S$ be an $\infty$-category.
First, we consider the $\infty$-category
${\rm Fun}(S,\oplaxmoncat{})$.
By Proposition~\ref{prop:description-fun-S-general-C},
we have an equivalence
${\rm Fun}(S,\oplaxmoncat{})\simeq 
\mathsf{Fun}(S,\oplaxmoncat{})$.
The inclusion functor
$\oplaxmoncat{}\hookrightarrow\cart{\Delta^{\natural}}$
induces a fully faithful functor
${\rm Fun}(S,\oplaxmoncat{})\to
{\rm Fun}(S,\cart{\Delta^{\natural}})$,
and we have an equivalence
${\rm Fun}(S,\cart{\Delta^{\natural}})\simeq
\mixed{(S^{\sharp},\Delta^{\natural})}$
by Proposition~\ref{prop:description-fun-S-cart-T-pair-version}.
%${\rm Fun}(S,\cart{\Delta_{\sharp}})\simeq
%\mathsf{Fun}(S,\cart{\Delta_{\sharp}})$
%of $\infty$-categories 
%by Proposition~\ref{prop:description-fun-S-general-C}.
Hence there is a fully faithful functor
\[ \mathsf{Fun}(S,\oplaxmoncat{})\longrightarrow
   \mixed{(S^{\sharp},\Delta^{\natural})}.\] 

We determine the essential image of
this functor.
Let
$p: X\to S\times \Delta$
be a mixed fibration over $(S,\Delta)$. 
%be an object of $\mixed{(S_{\sharp},\Delta_{\natural})}$.
%$\mathsf{Fun}(S,\cart{\Delta_{\sharp}})$.
%By Proposition~\ref{prop:cocart-cart-description},
%$p=\pi_S\circ P$ is a coCartesian fibration,
%$q=\pi_{\Delta}\circ P$ is a Cartesian fibration, 
%$P$ sends $p$-coCartesian morphisms to
%$\pi_S$-coCartesian morphisms, and
%$P$ sends $q$-Cartesian morphisms to
%$\pi_{\Delta}$-Cartesian morphisms.
%We denote by $X_{s,[m]}$
%the fiber of $P: X\to S\times\Delta$
%at $(s,[m])\in S\times \Delta$.
%We have a coCartesian fibration 
%$p_{[m]}: X_{[m]}\to S$ for any $[m]\in\Delta$.
%where $X_{\bullet,[m]}$
%is the fiber of $p_{\Delta}: X\to \Delta$ at $[m]\in\Delta$
%and $p_{\bullet,[m]}=p_{[m]}$.
%Since $p_{\Delta}: X\to \Delta$ is a Cartesian fibration
%and $p: X\to S\times \Delta$ preserves Cartesian
%morphisms over $\Delta$,
Any morphism $\phi: [m']\to [m]$ of $\Delta$
induces a functor 
$\phi^{!}: X_{[m]}\to X_{[m']}$ over $S$.
In particular,
we obtain a morphism
\[ X_{[m]}\longrightarrow
   \overbrace{X_{[1]}\times_{S}\cdots
              \times_{S}X_{[1]}}^m\]
in $\cocart{S^{\sharp}}$ for any $[m]\in \Delta$.
We also call it a Segal morphism. 
We consider the following condition on $p$\,:

\begin{condition}\rm\label{item:Segal-map-equivalence-over-S}
The Segal morphisms
$X_{[m]}\longrightarrow
   \overbrace{X_{[1]}\times_{S}\cdots
              \times_{S}X_{[1]}}^m$
are equivalences 
in the $\infty$-category $\cocart{S^{\sharp}}$
for all $[m]\in \Delta$.
\end{condition}

\begin{lemma}\label{lemma:mixed-fibration-S-segal-condition}
Let $p: X\to S\times \Delta$ be a mixed fibration
over $(S,\Delta)$.
If $p$ satisfies Condition~\ref{item:Segal-map-equivalence-over-S},
then $p$ is a marked mixed fibration over 
$(S^{\sharp},\Delta^{\natural})$.
\end{lemma}

\proof
If $p$ satisfies Condition~\ref{item:Segal-map-equivalence-over-S},
then the functor $(\rho^i)^!: X_{[m]}\to X_{[1]}$ 
over $S$ preserves
coCartesian morphisms for any $[m]\in \Delta$ and $1\le i\le m$.
For an inert morphism
$\phi: [m']\to [m]$ of $\Delta$,
there is a commutative diagram
\[ \begin{array}{ccc}
    X_{[m]}&\stackrel{\simeq}{\longrightarrow}&
    X_{[1]}\times_S\cdots\times_S X_{[1]}\\[1mm]
    \mbox{$\scriptstyle \phi^!$}\bigg\downarrow
    \phantom{\mbox{$\scriptstyle \phi^!$}}
    && \bigg\downarrow \\[3mm]
    X_{[m']}& \stackrel{\simeq}{\longrightarrow}&
    X_{[1]}\times_S\cdots\times_S X_{[1]},\\
   \end{array}\]
where the right vertical arrow is a projection.
This implies that $\phi^!: X_{[m]}\to X_{[m']}$
preserves coCartesian morphisms.
Hence $p$ is a marked mixed fibration
over $(S^{\sharp},\Delta^{\natural})$.
\qed

\begin{proposition}\label{prop:lands-in-oplaxmonoid-condition}
A mixed fibration $p: X\to S\times \Delta$
over $(S,\Delta)$
satisfies 
Condition~\ref{item:Segal-map-equivalence-over-S}
if and only if 
$p$ is an object
of $\mathsf{Fun}(S,\oplaxmoncat{})$.
%Let $p: X\to S\times \Delta$
%be a mixed fibration over $(S_{\sharp},\Delta_{\natural})$.
%%%be an object of $\mathfrak{Fun}(S,\cart{\Delta})$.
%If $p$ satisfies 
%Condition~\ref{item:Segal-map-equivalence-over-S},
%then $p$ is an object
%of $\mathsf{Fun}(S,\oplaxmoncat{})$.
%then $p_s: X_{s}\to \Delta$
%is an object of $\oplaxmoncat{}$
%for any $s\in S$.
\end{proposition}

\proof
%Note that Condition~\ref{item:Segal-map-equivalence-over-S}
%implies that the functor 
%$\phi^{!}: X_{[m']}\to X_{[m]}$
%sends $p_{[m']}$-coCartesian morphisms to
%$p_{[m]}$-coCartesian morphisms
%if $\phi: [m]\to [m']$ is an inert morphism
%of $\Delta$. 

%We have a Cartesian fibration $p_s: X_{s}\to \Delta$
%for any $s\in S$.
%where $X_{s,\bullet}$ is the fiber of $p: X\to S$ at $s\in S$.
First, we suppose that $p$ satisfies 
Condition~\ref{item:Segal-map-equivalence-over-S}.
Then the Segal morphism
$X_{(s,[m])}\to X_{(s,[1])}\times\cdots \times X_{(s,[1])}$
is an equivalence of $\infty$-categories
for any $(s,[m])\in S\times \Delta$.
Furthermore,
since $p: X\to S\times\Delta$
is a marked mixed fibration over 
$(S^{\sharp},\Delta^{\natural})$
by Lemma~\ref{lemma:mixed-fibration-S-segal-condition},
the functor $e_{!}: X_{s}\to X_{s'}$ over $\Delta$
induced by any morphism $e: s\to s'$ of $S$
carries $p_s$-Cartesian morphisms
over inert morphisms of $\Delta$ 
to $p_{s'}$-Cartesian morphisms.
Thus, $p$ is an object of 
$\mathsf{Fun}(S,\oplaxmoncat{})$.

Conversely,
we suppose that $p$ is an object of 
$\mathsf{Fun}(S,\oplaxmoncat{})$.
Since $p: X\to S\times \Delta$ is a marked mixed fibration
over $(S^{\sharp},\Delta^{\natural})$,
the functor $\phi^!: X_{[m]}\to X_{[m']}$
over $S$ induced by any inert morphism
$\phi: [m']\to [m]$ of $\Delta$
preserves coCartesian morphisms.
Hence the Segal morphism
$X_{[m]}\to X_{[1]}\times_S\cdots\times_S X_{[1]}$
preserves coCartesian morphisms over $S$
for any $[m]\in\Delta$.
Furthermore,
the maps 
$X_{(s,[m])}\to X_{(s,[1])}\times\cdots X_{(s,[1])}$
on fibers 
are equivalences by the assumption
for all $s\in S$.
By the dual of \cite[Proposition~3.3.1.5]{Lurie1},
we see that $p$ satisfies 
Condition~\ref{item:Segal-map-equivalence-over-S}.
\qed

%Since $p: X\to S$ is a coCartesian fibration
%and $p: X\to S\times \Delta$ preserves coCartesian
%morphisms over $S$,
%any morphism $e: s\to s'$ of $S$
%induces a functor 
%$e_{!}: X_{s}\to X_{s'}$
%over $\Delta$.
%Since $p: X\to S\times\Delta$
%is a mixed fibration over $(S_{\sharp},\Delta_{\natural})$,
%the functor $e_{!}: X_{s}\to X_{s'}$
%induced by a morphism $e: s\to s'$ of $S$
%carries $p_s$-Cartesian morphisms
%over inert morphisms of $\Delta$ 
%to $p_{s'}$-Cartesian morphisms.

%\begin{theorem}
\begin{corollary}\label{cor:description-fun-S-oplaxmoncat}
The $\infty$-category
${\rm Fun}(S,\oplaxmoncat)$
is equivalent to a full subcategory
of $\mixed{(S_{\sharp},\Delta_{\natural})}$
spanned by marked mixed fibrations
%over $(S_{\sharp},\Delta_{\natural})$
over $(S^{\sharp},\Delta^{\natural})$
which satisfy Condition~\ref{item:Segal-map-equivalence-over-S}.
%The $\infty$-category $\mathfrak{Fun}(S,\oplaxmoncat{})$
%is given as follows. 
%The objects of $\mathfrak{Fun}(S,\oplaxmoncat{})$
%are functors $P: X\to S\times\Delta$ of $\infty$-categories
%such that $P$ satisfies 
%Conditions~\ref{item:coCart-preserve-coCart},
%\ref{item:Cart-preserve-Cart}, 
%and \ref{item:Segal-map-equivalence-over-S}.
%The morphisms between 
%objects $P: X\to S\times\Delta$
%and $Q: X\to S\times \Delta$ of $\mathfrak{Fun}(S,\oplaxmoncat{})$
%are functors $F: P\to Q$ of $\mathfrak{Fun}(S,\cart{\Delta})$
%such that $F_s: X_s\to Y_s$ preserves
%Cartesian morphisms over inert morphisms of $\Delta$
%for any $s\in S$.
\end{corollary}
%\end{theorem}

Next, we consider 
a description of the $\infty$-category
${\rm Mon}(\oplaxmoncat{})$
in terms of mixed fibrations.
By definition,
%the $\infty$-category 
${\rm Mon}(\oplaxmoncat{})$
is a full subcategory 
of ${\rm Fun}(\Delta^{\rm op},\oplaxmoncat{})$,
and hence 
there is a fully faithful functor
\[ {\rm Mon}(\oplaxmoncat{})\longrightarrow
   \mixed{((\Delta^{\rm op})^{\sharp},(\Delta)^{\natural})} \]
by Corollary~\ref{cor:description-fun-S-oplaxmoncat}.
%spanned by mixed fibrations over 
%$(\Delta^{\rm op}_{\sharp},\Delta_{\natural})$
%which satisfies 
%Condition~\ref{item:Segal-map-equivalence-over-S}
%by Corollary~\ref{cor:description-fun-S-oplaxmoncat}.
We shall determine the essential image of this functor.

Let 
%$p: X\longrightarrow \Delta_1^{\rm op}\times\Delta_{\otimes}$
$p: X\to \Delta_{\boxtimes}^{\rm op}\times\Delta_{\otimes}$
be a mixed fibration over $(\Delta_{\boxtimes}^{\rm op},\Delta_{\otimes})$,
%\[ p: X\longrightarrow \Delta_1^{\rm op}\times\Delta_{\otimes} \]
%be a mixed fibration over $(\Delta_1^{\rm op},\Delta_{\otimes})$,
where $\Delta_{\boxtimes}=\Delta_{\otimes}=\Delta$.
%a functor of $\infty$-categories,
%which satisfies Conditions~\ref{item:coCart-preserve-coCart}
%and \ref{item:Cart-preserve-Cart},
%where $\Delta_1=\Delta_{\otimes}=\Delta$. 
We denote by $X_{[m],[n]}$
%the fiber of $p$ at $([m],[n])\in \Delta^{\rm op}\times\Delta$.
the fiber of $p$ at $([m],[n])\in 
\Delta_{\boxtimes}^{\rm op}\times\Delta_{\otimes}$.
Recall that we have coCartesian fibrations
%$p_{\bullet,[m]}:X_{\bullet,[m]}\to\Delta^{\rm op}$
$p_{\bullet,[n]}:X_{\bullet,[n]}\to\Delta_{\boxtimes}^{\rm op}$
%for $[n]\in\Delta$,
for each $[n]\in\Delta_{\otimes}$,
and 
Cartesian fibrations
%$p_{[n],\bullet}: X_{[n],\bullet}\to\Delta$
$p_{[m],\bullet}: X_{[m],\bullet}\to\Delta_{\otimes}$
%for $[m]\in\Delta^{\rm op}$.
for each $[m]\in\Delta_{\boxtimes}^{\rm op}$.
We consider the following conditions on $p$\,:

%\begin{enumerate}
%\item 
%The composite $\pi_1\circ p: X\to \Delta_1^{\rm op}$
%is a coCartesian fibration,
%where $\pi_1:\Delta_1^{\rm op}\times\Delta_{\otimes}\to\Delta_1^{\rm op}$
%is the projection. 
%When we regard $\pi_1$ as a coCartesian fibration,
%the functor $p$ preserves coCartesian morphisms.

%\item
%The composite $\pi_2\circ p: X\to \Delta_{\otimes}$
%is a Cartesian fibration,
%where $\pi_2:\Delta_1^{\rm op}\times\Delta_{\otimes}\to\Delta_{\otimes}$
%is the projection. 
%When we regard $\pi_2$ as a Cartesian fibration,
%the functor $p$ preserves Cartesian morphisms.

\begin{condition}\rm\label{item:coCart-Segal-map-again}
The Segal morphisms
%\[ X_{\bullet,[m]}\longrightarrow
%   \overbrace{X_{\bullet,[1]}\times_{\Delta^{\rm op}}\cdots
%              \times_{\Delta^{\rm op}}X_{\bullet,[1]}}^m\]
\[ X_{\bullet,[n]}\longrightarrow
   \overbrace{X_{\bullet,[1]}\times_{\Delta_{\boxtimes}^{\rm op}}\cdots
              \times_{\Delta_{\boxtimes}^{\rm op}}X_{\bullet,[1]}}^n\]
are equivalences in the $\infty$-category
%$\cocart{\Delta^{\rm op}}$
$\cocart{(\Delta_{\boxtimes}^{\rm op})^{\sharp}}$
%of coCartesian fibrations
%over $\Delta^{\rm op}$
%for all $[m]\in\Delta$.
for all $[n]\in\Delta_{\otimes}$.
\end{condition}

\begin{condition}\rm\label{item:Cart-Segal-map-again}
The Segal morphisms
%\[ X_{[m],\bullet}\longrightarrow
%   \overbrace{X_{[1],\bullet}\times_{\Delta}\cdots
%              \times_{\Delta}X_{[1],\bullet}}^n\]
\[ X_{[m],\bullet}\longrightarrow
   \overbrace{X_{[1],\bullet}\times_{\Delta_{\otimes}}\cdots
              \times_{\Delta_{\otimes}}X_{[1],\bullet}}^m\]
are equivalences in the $\infty$-category
%$\cart{\Delta}$
$\cart{(\Delta_{\otimes})^{\sharp}}$
%of Cartesian fibrations
%over $\Delta$
%for all $[n]\in\Delta^{\rm op}$.
for all $[m]\in\Delta_{\boxtimes}^{\rm op}$.
\end{condition}

%Let $P: X\to \Delta_1^{\rm op}\times\Delta_{\otimes}$
%be a functor of $\infty$-categories
%satisfying Conditions~\ref{item:coCart-preserve-coCart}
%and \ref{item:Cart-preserve-Cart}.
%We regard $P$ as an object of
%$\mathfrak{Fun}(\Delta_1^{\rm op},\cart{\Delta_{\otimes}})$.

\begin{proposition}
%Let $p: X\to\Delta^{\rm op}\times\Delta$
Let $p: X\to\Delta_{\boxtimes}^{\rm op}\times\Delta_{\otimes}$
be a mixed fibration over 
%$(\Delta^{\rm op},\Delta)$.
$(\Delta_{\boxtimes}^{\rm op},\Delta_{\otimes})$.
%Under the equivalence
%$\mathsf{Fun}(\Delta^{\rm op},\cart{\Delta_{\natural}})\simeq
%{\rm Fun}(\Delta^{\rm op},\cart{\Delta_{\natural}})$
%of Proposition~\ref{prop:description-fun-S-general-C},
The map $p$ is in the essential image of the
functor 
%${\rm Mon}(\oplaxmoncat{})\to
%   \mixed{(\Delta_{\sharp}^{\rm op},\Delta_{\natural})}$
${\rm Mon}(\oplaxmoncat{})\to
   \mixed{((\Delta_{\boxtimes}^{\rm op})^{\sharp},
           (\Delta_{\otimes})^{\natural})}$
if and only if
$p$ satisfies Conditions~\ref{item:coCart-Segal-map-again} and
\ref{item:Cart-Segal-map-again}.
\end{proposition}

\proof
First, we suppose that $p$ is in the essential image 
of the functor 
${\rm Fun}(\Delta_{\boxtimes}^{\rm op},\oplaxmoncat)\to
\mixed{((\Delta_{\boxtimes}^{\rm op})^{\sharp},
(\Delta_{\otimes})^{\natural})}$.
By Corollary~\ref{cor:description-fun-S-oplaxmoncat},
$p$ satisfies Condition~\ref{item:coCart-Segal-map-again}.
By the assumption,
the Segal morphisms
$X_{[m],\bullet}\to
   X_{[1],\bullet}\times_{\Delta_{\otimes}}\cdots
              \times_{\Delta_{\otimes}}X_{[1],\bullet}$
are equivalences in $\oplaxmoncat{}$ 
for all $[m]\in\Delta^{\rm op}_{\boxtimes}$.
This implies that Condition~\ref{item:Cart-Segal-map-again} 
holds.

Conversely,
we suppose that $p$ satisfies 
Conditions~\ref{item:coCart-Segal-map-again} and
\ref{item:Cart-Segal-map-again}.
By Corollary~\ref{cor:description-fun-S-oplaxmoncat},
%Proposition~\ref{prop:lands-in-oplaxmonoid-condition},
%and
%Lemma~\ref{lemma:preserves-inert-morphisms-condition},
%the functor $\Delta_1^{\rm op}\to \cart{\Delta_{\otimes}}$
%corresponding to 
$p$ is in the essential image 
of the functor
${\rm Fun}(\Delta_{\boxtimes}^{\rm op},\oplaxmoncat)\to
\mixed{((\Delta_{\boxtimes}^{\rm op})^{\sharp},
(\Delta_{\otimes})^{\natural})}$.
Since %the inclusion functor
%$\oplaxmoncat{}\hookrightarrow \cart{\Delta_{\otimes}}$ 
%is conserveative,
$\oplaxmoncat{}$ is a full subcategory
of $\cart{(\Delta_{\otimes})^{\natural}}$,
Condition~\ref{item:Cart-Segal-map-again}
implies that the Segal morphisms
%$X_{[n],\bullet}\to
%   X_{[1],\bullet}\times_{\Delta^{\rm op}}\cdots
%              \times_{\Delta^{\rm op}}X_{[1],\bullet}$
$X_{[m],\bullet}\to
   X_{[1],\bullet}\times_{\Delta_{\otimes}}\cdots
              \times_{\Delta_{\otimes}}X_{[1],\bullet}$
are equivalences in $\oplaxmoncat{}$
%$[n]\in\Delta^{\rm op}$. 
for all $[m]\in\Delta_{\boxtimes}^{\rm op}$. 
Hence $p$ is in the essential image of the functor
${\rm Mon}(\Delta_{\boxtimes}^{\rm op},\oplaxmoncat)\to
\mixed{((\Delta_{\boxtimes}^{\rm op})^{\sharp},(\Delta_{\otimes})^{\natural})}$.
%Thus, $p$ corresponds to a duoidal $\infty$-category.
\qed

\begin{definition}\rm
We also say that
a mixed fibration
\[ p: X\to \Delta_{\boxtimes}^{\rm op}\times\Delta_{\otimes} \]
over $(\Delta_{\boxtimes}^{\rm op},\Delta_{\otimes})$
is a duoidal $\infty$-category 
if it satisfies Conditions~\ref{item:coCart-Segal-map-again} 
and \ref{item:Cart-Segal-map-again}.
\end{definition}

\begin{example}\rm
If $R$ is a commutative ring,
then the category of $R$-$R$-bimodules has
the structure of a duoidal category.
We consider a generalization of this fact
in $\infty$-categories.
Let $\mathcal{C}$ be a presentable
$\mathbb{E}_2$-monoidal $\infty$-category,
and let $A$ be an $\mathbb{E}_2$-algebra object of $\mathcal{C}$.
There is a duoidal $\infty$-category
$p: X\to \Delta_{\boxtimes}^{\rm op}\times
\Delta_{\otimes}$ in which
$X_{[1],[1]}$ is equivalent to
the $\infty$-category ${}_A{\rm BMod}_A(\mathcal{C})$
of $A$-$A$-bimodules in $\mathcal{C}$.
For $M,N\in {}_A{\rm BMod}_A(\mathcal{C})$,
we have $M\otimes N\simeq M\otimes_AN$
and $M\boxtimes N\simeq A\otimes_{A\otimes A}(M\otimes N)
\otimes_{A\otimes A}A$.
See \cite{Torii4} for more details
and a generalization to operadic modules.
\end{example}

%Recall that the functor $R: \cat\to\cat$
%assigns to an $\infty$-category
%its opposite $\infty$-category. 
We shall show that the duoidal $\infty$-categories
are closed under the functor $(-)^{\rm op}$.

\begin{proposition}
Let $p: X\to \Delta_{\boxtimes}^{\rm op}\times\Delta_{\otimes}$
be a mixed fibration over $(\Delta_{\boxtimes}^{\rm op},\Delta_{\otimes})$.
Then $p$ is a duoidal $\infty$-category
if and only if 
$p^{\rm op}: X^{\rm op}\to \Delta_{\otimes}^{\rm op}\times\Delta_{\boxtimes}$
is a duoidal $\infty$-category. 
\end{proposition}

\proof
By Remark~\ref{remark:mixed-fib-R-duality},
%By Remark~\ref{remark:cocart-cart-duality},
if $p$ is a mixed fibration over
$(\Delta_{\boxtimes}^{\rm op},\Delta_{\otimes})$,
%$\mathfrak{Fun}(\Delta_{\boxtimes}^{\rm op},\cart{\Delta_{\otimes}})$
%if and only if $R(p)$ is an object of
then $p^{\rm op}$ is a mixed fibration over
$(\Delta_{\otimes}^{\rm op},\Delta_{\boxtimes})$.
%$\mathfrak{Fun}(\Delta_{\otimes}^{\rm op},\cart{\Delta_{\boxtimes}})$.
By the symmetry of Conditions~\ref{item:coCart-Segal-map-again}
and \ref{item:Cart-Segal-map-again},
we see that $p$ is a duoidal $\infty$-category
if and only if $p^{\rm op}$ is a duoidal $\infty$-category.
\qed

%\newpage

\subsection{Bilax monoidal functors}
\label{subsection:bilax-monoidal-functor}

Let $p:X\to \Delta_{\boxtimes}^{\rm op}\times\Delta_{\otimes}$ and 
$q: Y\to \Delta_{\boxtimes}^{\rm op}\times\Delta_{\otimes}$
be duoidal $\infty$-categories.
Suppose we have a functor $f:X\to Y$ 
over $\Delta_{\boxtimes}^{\rm op}\times\Delta_{\otimes}$
which is a morphism of ${\rm Mon}(\oplaxmoncat{})$.
%$p:X\to \Delta_{\boxtimes}^{\rm op}\times\Delta_{\otimes}$ and 
%$q: Y\to \Delta_{\boxtimes}^{\rm op}\times\Delta_{\otimes}$.
Then we obtain strong monoidal functors
$f_{\bullet,[n]}: X_{\bullet,[n]}\to Y_{\bullet,[n]}$
for each $[n]\in\Delta_{\otimes}$,
and oplax monoidal functors
$f_{[m],\bullet}: X_{[m],\bullet}\to Y_{[m],\bullet}$
for each $[m]\in\Delta_{\boxtimes}^{\rm op}$.
This is asymmetric under $(-)^{\rm op}$.
We shall introduce a bilax monoidal 
functor between duoidal $\infty$-categories,
which is symmetric under $(-)^{\rm op}$. 

\begin{definition}\rm
Let 
$p: X\to \Delta_{\boxtimes}^{\rm op}\times \Delta_{\otimes}$
and $q: Y\to \Delta_{\boxtimes}^{\rm op}\times \Delta_{\otimes}$
be duoidal $\infty$-categories,
and let $f: X\to Y$ be a functor
over $\Delta_{\boxtimes}^{\rm op}\times\Delta_{\boxtimes}$.
%Let 
%$p: X\to \Delta_{\boxtimes}^{\rm op}\times \Delta_{\otimes}$
%and $q: Y\to \Delta_{\boxtimes}^{\rm op}\times \Delta_{\otimes}$
%be duoidal $\infty$-categories, and
%let 
%$\pi_1: \Delta_{\boxtimes}^{\rm op}\times\Delta_{\otimes}\to \Delta_{\boxtimes}^{\rm op}$
%and $\pi_2: \Delta_{\boxtimes}^{\rm op}\times\Delta_{\otimes}\to\Delta_{\otimes}$
%be the projections.
We say that $f$ is a bilax monoidal functor 
between $p$ and $q$ 
if $f$ is a morphism
of $\mixed{((\Delta_{\boxtimes}^{\rm op})^{\natural},
(\Delta_{\otimes})^{\natural})}$.
\end{definition}

By definition, 
a bilax monoidal functor $f:X\to Y$
between duoidal $\infty$-categories 
$p: X\to\Delta_{\boxtimes}^{\rm op}\times \Delta_{\otimes}$ and 
$q: Y\to\Delta_{\boxtimes}^{\rm op}\times\Delta_{\otimes}$
induces lax monoidal functors
\[ f_{\bullet,[n]}: X_{\bullet,[n]}\longrightarrow
                  Y_{\bullet,[n]} \]
for each $[n]\in\Delta_{\otimes}$,
and oplax monoidal functors
\[ f_{[m],\bullet}: X_{[m],\bullet}\longrightarrow
                  Y_{[m],\bullet} \]
for each $[m]\in\Delta_{\boxtimes}^{\rm op}$.  
By the symmetry of the definition of bilax monoidal
functors,
we easily obtain the following proposition.

\begin{proposition}
If $f$ is a bilax monoidal functor
between duoidal $\infty$-categories,
then $f^{\rm op}$ is also
a bilax monoidal functor.
\end{proposition}

\begin{definition}\rm
We define 
\[ \duo^{\rm bilax} \]
to be the full subcategory
of $\mixed{((\Delta_{\boxtimes}^{\rm op})^{\natural},
(\Delta_{\otimes})^{\natural})}$
spanned by duoidal $\infty$-categories.
We call $\duo^{\rm bilax}$
the $\infty$-category of 
duoidal $\infty$-categories and bilax monoidal functors.
\end{definition}

\begin{remark}\rm
The functor $(-)^{\rm op}: \cat\to\cat$
induces an equivalence 
\[ (-)^{\rm op}: \duo^{\rm bilax}\stackrel{\simeq }
      {\longrightarrow}\duo^{\rm bilax} \]
of $\infty$-categories.
\end{remark}

For a duoidal $\infty$-category
$p: X\to\Delta_{\boxtimes}^{\rm op}\times \Delta_{\otimes}$,
we call the fiber $X_{[1],[1]}$
the underlying $\infty$-category.
%the duoidal $\infty$-category $p$.
Assigning to a duoidal $\infty$-category
its underlying $\infty$-category,
we obtain a functor
\[ {\rm ev}_{[1],[1]}: \duo^{\rm bilax}\longrightarrow \cat \]  
of $\infty$-categories.

We shall show that
%the $\infty$-category 
$\duo^{\rm bilax}$
%of duoidal $\infty$-categories 
is equivalent to 
the $\infty$-category
$\laxmon(\oplaxmoncat{})$
of monoid objects and lax monoidal functors
in $\oplaxmoncat{}$.

\begin{theorem}\label{theorem:duo-laxmon-oplaxmon}
There is an equivalence 
\[ \duo^{\rm bilax} \simeq \laxmon(\oplaxmoncat{})\]
of $\infty$-categories.
\end{theorem}

\proof
The objects of $\duo^{\rm bilax}$
coincide with those of
$\laxmon(\oplaxmoncat{})$.
We shall show that 
the morphisms in $\duo^{\rm bilax}$ coincide with
those in $\laxmon(\oplaxmoncat{})$.

Let $p: X\to\Delta_{\boxtimes}^{\rm op}\times\Delta_{\otimes}$
and $q: Y\to\Delta_{\boxtimes}^{\rm op}\times\Delta_{\otimes}$
be duoidal $\infty$-categories. 
The morphisms between $p$ and $q$
in $\duo^{\rm bilax}$ are functors
$f: X\to Y$ over $\Delta_{\boxtimes}^{\rm op}\times \Delta_{\otimes}$
such that 
$f$ carries $p_{\Delta_{\boxtimes}^{\rm op}}$-coCartesian morphisms
over inert morphisms of $\Delta^{\rm op}_{\boxtimes}$
to $q_{\Delta_{\boxtimes}^{\rm op}}$-coCartesian morphisms,
and that $p_{\Delta_{\otimes}}$-Cartesian morphisms
over inert morphisms of $\Delta_{\otimes}$
to $q_{\Delta_{\otimes}}$-Cartesian morphisms.
On the other hand,
the morphisms between $p$ and $q$
in $\laxmon(\oplaxmoncat{})$
are functors 
$f: X\to Y$ over $\Delta_{\boxtimes}^{\rm op}\times \Delta_{\otimes}$
such that 
$f$ carries $p_{\Delta_{\boxtimes}^{\rm op}}$-coCartesian morphisms
over inert morphisms of $\Delta^{\rm op}_{\boxtimes}$
to $q_{\Delta_{\boxtimes}^{\rm op}}$-coCartesian morphisms,
and that the functor $f_{[m],\bullet}: X_{[m],\bullet}\to 
Y_{[m],\bullet}$ over $\Delta_{\otimes}$
is a morphism of $\oplaxmoncat{}$
for each $[m]\in\Delta^{\rm op}_{\boxtimes}$
by Definition~\ref{def:lax-mon-cat-C}.
Recall that
$f_{[m],\bullet}$ is a morphism of $\oplaxmoncat{}$
if and only if $f_{[m],\bullet}$
sends $p_{[m],\bullet}$-Cartesian morphisms 
over inert morphisms of $\Delta_{\otimes}$
to $q_{[m],\bullet}$-Cartesian morphisms.
By Remark~\ref{remark:p-cart-equivalent-pt-cart},
the morphisms between $p$ and $q$ in $\laxmon(\oplaxmoncat{})$
coincide with
those in $\duo^{\rm bilax}$.
This completes the proof.
\qed

\begin{remark}\rm
By the symmetry of the definition
of duoidal $\infty$-categories and bilax monoidal functors,
there is an equivalence 
\[ \duo^{\rm bilax}\simeq 
   \mathsf{Mon}^{\rm oplax}(\mathsf{Mon}^{\rm lax}(\cat)) \]
of $\infty$-categories. 
\end{remark}

%\newpage

\subsection{Double lax and oplax monoidal functors}

In \S\ref{subsection:bilax-monoidal-functor}
we introduced bilax monoidal functors
between duoidal $\infty$-categories.
We have two other kinds of functors
between duoidal $\infty$-categories.
We will call these functors double lax monoidal functors
and double oplax monoidal functors.
In order to formulate these functors,
we shall give other descriptions
of duoidal $\infty$-categories.
Using one of these descriptions,
we will define an $\infty$-category
of duoidal $\infty$-categories and double lax monoidal functors.
By duality,
we will also define an $\infty$-category 
of duoidal $\infty$-categories and double oplax monoidal functors.

By Definition~\ref{definition:duoidal-infty-category-monoid-oplax},
a duoidal $\infty$-category is a 
monoid object of $\oplaxmoncat$. 
Since there is an equivalence
$(-)^{\rm op}: \oplaxmoncat\stackrel{\sim}{\to}\laxmoncat$
of $\infty$-categories,
we can regard a duoidal $\infty$-category
as a monoid object of $\laxmoncat$.
Furthermore,
by Proposition~\ref{prop:moncal-monsf},
we can identify a duoidal $\infty$-category
with an object of $\mathsf{Mon}(\laxmoncat)$,
which is a categorical fibration
$p: X\to \Delta_{\otimes}^{\rm op}\times\Delta_{\boxtimes}^{\rm op}$
satisfying the following conditions:

\begin{condition}\label{condition:duoidal-monoid-lax}\rm
Let
$p: X\to \Delta_{\otimes}^{\rm op}\times\Delta_{\boxtimes}^{\rm op}$
be a categorical fibration.
\begin{enumerate}
\item The composite $\pi_{\otimes}\circ p: 
X\to \Delta_{\otimes}^{\rm op}$
is a coCartesian fibration,
and $p$ carries $(\pi_{\otimes}\circ p)$-coCartesian 
morphisms to $\pi_{\otimes}$-coCartesian morphisms.
\item
For each $[n]\in \Delta_{\otimes}^{\rm op}$,
the restriction
$p_{[n],\bullet}: X_{[n],\bullet}\to \Delta_{\boxtimes}^{\rm op}$
is a monoidal $\infty$-category,
\item
For each morphism
$[n]\to [n']$ in $\Delta_{\otimes}^{\rm op}$,
the induced map
$X_{[n],\bullet}\to X_{[n'],\bullet}$
over $\Delta_{\boxtimes}^{\rm op}$
is a lax monoidal functor.
\item
For each $[n]\in\Delta_{\otimes}^{\rm op}$,
the Segal morphism
\[ X_{[n],\bullet}\longrightarrow
   \overbrace{X_{[1],\bullet}\times_{\Delta_{\boxtimes}^{\rm op}}\cdots
              \times_{\Delta_{\boxtimes}^{\rm op}}X_{[1],\bullet}}^{n} \]
is an equivalence of coCartesian fibrations
over $\Delta_{\boxtimes}^{\rm op}$.
\end{enumerate}
\end{condition}

\begin{definition}\rm
We also say that a categorical fibration 
\[ p: X\to \Delta_{\otimes}^{\rm op}\times\Delta_{\boxtimes}^{\rm op} \]
is a duoidal $\infty$-category
if it satisfies 
Condition~\ref{condition:duoidal-monoid-lax}.
\end{definition}

\begin{definition}\rm
Let $p: X\to \Delta_{\otimes}^{\rm op}
\times\Delta_{\boxtimes}^{\rm op}$
and 
$q: Y\to \Delta_{\otimes}^{\rm op}
\times\Delta_{\boxtimes}^{\rm op}$
be duoidal $\infty$-categories.
Let $f: X\to Y$ be a functor over
$\Delta_{\otimes}^{\rm op}
\times\Delta_{\boxtimes}^{\rm op}$.
We say that $f$ is a double lax monoidal functor
if it carries $(\pi_{\otimes}\circ p)$-coCartesian morphisms
over inert morphisms of $\Delta_{\otimes}^{\rm op}$
to $(\pi_{\otimes}\circ q)$-coCartesian morphisms
and if the induced functor
$f_{[n],\bullet}: X_{[n],\bullet}\to Y_{[n],\bullet}$
over $\Delta_{\boxtimes}^{\rm op}$
is a lax monoidal functor for each $[n]\in\Delta_{\otimes}^{\rm op}$.
\end{definition}

In other words,
a double lax monoidal functor is a morphism
in $\laxmon(\laxmoncat)$.

\begin{definition}\rm
We define an $\infty$-category
\[ \duo^{\rm dlax} \]
to be the $\infty$-category
$\laxmon(\laxmoncat)$.
We call $\duo^{\rm dlax}$ 
the $\infty$-category of duoidal $\infty$-categories
and double lax monoidal functors.
\end{definition}

Next,
by duality,
we introduce double oplax monoidal functors
between duoidal $\infty$-categories.

\begin{definition}\rm
We define an $\infty$-category
\[ \duo^{\rm doplax} \]
to be the $\infty$-category
$\oplaxmon(\oplaxmoncat)$.
We also say that 
an object of $\duo^{\rm doplax}$ is a duoidal $\infty$-category.
A morphism
in $\duo^{\rm doplax}$ is said to be a double 
oplax monoidal functor.
We call $\duo^{\rm doplax}$ 
the $\infty$-category of duoidal $\infty$-categories
and double oplax monoidal functors.
\end{definition}

\begin{remark}\rm
\label{remark:dual-Cartesian-description-duoidal-category}
By Definition~\ref{def:description-oplaxmoncat},
an object of $\duo^{\rm doplax}$
is a categorical fibration
$p: X\to \Delta_{\otimes}\times\Delta_{\boxtimes}$
which satisfies the following conditions:

\begin{enumerate}
\item The composite $\pi_{\otimes}\circ p: 
X\to \Delta_{\otimes}$
is a Cartesian fibration,
and $p$ carries $(\pi_{\otimes}\circ p)$-Cartesian 
morphisms to $\pi_{\otimes}$-Cartesian morphisms.
\item
For each $[n]\in \Delta_{\otimes}$,
the restriction
$p_{[n],\bullet}: X_{[n],\bullet}\to \Delta_{\boxtimes}$
is a monoidal $\infty$-category.
\item
For each morphism
$[n']\to [n]$ in $\Delta_{\otimes}$,
the induced map
$X_{[n],\bullet}\to X_{[n'],\bullet}$
over $\Delta_{\boxtimes}$
is an oplax monoidal functor.
\item
For each $[n]\in\Delta_{\otimes}$,
the Segal morphism
\[ X_{[n],\bullet}\longrightarrow
   \overbrace{X_{[1],\bullet}\times_{\Delta_{\boxtimes}}\cdots
              \times_{\Delta_{\boxtimes}}X_{[1],\bullet}}^{n} \]
is an equivalence of Cartesian fibrations
over $\Delta_{\boxtimes}$.
\end{enumerate}
%\end{condition}

%\begin{definition}\rm
Let $p: X\to \Delta_{\otimes}\times\Delta_{\boxtimes}$
and 
$q: Y\to \Delta_{\otimes}\times\Delta_{\boxtimes}$
be duoidal $\infty$-categories.
A double oplax monoidal functor from $p$ to $q$
is a functor $f: X\to Y$ over
$\Delta_{\otimes}\times\Delta_{\boxtimes}$
such that it carries $\pi_{\otimes}\circ p$-Cartesian morphisms
over inert morphisms of $\Delta_{\otimes}$
to $\pi_{\otimes}\circ q$-Cartesian morphisms
and that the induce functor
$f_{[n],\bullet}: X_{[n],\bullet}\to Y_{[n],\bullet}$
over $\Delta_{\boxtimes}$
is an oplax monoidal functor for each $[n]\in\Delta_{\otimes}$.
%\end{definition}
\end{remark}

\begin{remark}\rm
For a duoidal $\infty$-category
$p: X\to \Delta_{\otimes}^{\rm op}\times
\Delta_{\boxtimes}^{\rm op}$,
the map $p^{\rm op}: X^{\rm op}\to
\Delta_{\otimes}\times\Delta_{\boxtimes}$
is also a duoidal $\infty$-category.
This induces an equivalence
\[ (-)^{\rm op}: \duo^{\rm dlax}\stackrel{\simeq}{\longrightarrow}
      \duo^{\rm doplax} \]
%\[ R: \laxmon(\laxmoncat)\stackrel{\simeq}{\longrightarrow}
%      \oplaxmon(\oplaxmoncat) \]
of $\infty$-categories.
\end{remark}

\section{Bimonoids, double monoids, and
double comonoids in duoidal $\infty$-categories}
\label{section:bialgebra}

We can consider bimonoids,
double monoids, and double comonoids
in duoidal categories.
In this section we introduce notions of bimonoids,
double monoids, and double comonoids
in duoidal $\infty$-categories.
For a duoidal $\infty$-category 
$p: X\to \Delta_{\boxtimes}^{\rm op}\times\Delta_{\otimes}$,
we construct a monoidal structure
on the $\infty$-category $\alg{}{}^{\boxtimes}(X)$
of algebra objects of $(X,\boxtimes, 1_{\boxtimes})$
by using the monoidal structure $(\otimes, 1_{\otimes})$.
Similarly,
we construct a monoidal structure  
on the $\infty$-category $\coalg{}{}^{\otimes}(X)$
of coalgebra objects of $(X,\otimes, 1_{\otimes})$
by using the monoidal structure 
$(\boxtimes, 1_{\boxtimes})$.
We show that 
bimonoid objects of $X$ are equivalent to 
algebra objects of $\coalg{}{}^{\otimes}(X)$,
and that they are also
equivalent to coalgebra objects of $\alg{}{}^{\boxtimes}(X)$
by duality.
We also prove similar results on
double monoids and double comonoids
in duoidal $\infty$-categories.

\subsection{Algebra and coalgebra objects
in monoidal $\infty$-categories}
\label{subsection:algebra-coalgebra-monoida-category}

In this subsection
we recall algebra and coalgebra objects
in monoidal $\infty$-categories
(cf.~\cite[\S2.1.2]{Lurie2} and \cite[\S4.3]{Groth}).

First, we recall algebra objects
in a monoidal $\infty$-category.
Let $p: X\to \Delta^{\rm op}$ be a monoidal
$\infty$-category
with underlying $\infty$-category $X_{[1]}$.
%For simplicity,
%we say that $X_{[1]}$ is a monoidal $\infty$-category.
An algebra object of the monoidal $\infty$-category $X_{[1]}$ is a 
section $s: \Delta^{\rm op}\to X$ of $p$
such that $s$ takes inert morphisms of $\Delta^{\rm op}$
to $p$-coCartesian morphisms of $X$.

We can describe algebra objects
of $X_{[1]}$ in another way.
Let ${\rm id}_{\Delta^{\rm op}}: \Delta^{\rm op}\to\Delta^{\rm op}$
be the identity functor of $\Delta^{\rm op}$.
Note that ${\rm id}_{\Delta^{\rm op}}$ is a final object
of $\laxmoncat{}$.
An algebra object of $X_{[1]}$
is a morphism from ${\rm id}_{\Delta^{\rm op}}$ to $p$
in $\laxmoncat{}$.

Let ${\rm Fun}_{\Delta^{\rm op}}(\Delta^{\rm op},X)$
be the $\infty$-category of sections of $p$.
We denote by
\[ {\rm Alg}(X) \]
%\[ {\rm Alg}(X_{[1]}) \]
the full subcategory of ${\rm Fun}_{\Delta^{\rm op}}(\Delta^{\rm op},X)$
spanned by algebra objects of $X_{[1]}$.
Let $f$ be a lax monoidal functor from
$p: X\to\Delta^{\rm op}$ to $q: Y\to\Delta^{\rm op}$.
Since $f$ carries $p$-coCartesian morphisms
over inert morphisms of $\Delta^{\rm op}$
to $q$-coCartesian morphisms,
it induces a functor
\[ {\rm Alg}(f): 
   {\rm Alg}(X)\longrightarrow {\rm Alg}(Y) \]
%\[ {\rm Alg}(F_{[1]}): 
%   {\rm Alg}(X_{[1]})\longrightarrow {\rm Alg}(Y_{[1]}) \]
of $\infty$-categories.
Assigning 
to a monoidal $\infty$-category $p: X\to\Delta^{\rm op}$   
the $\infty$-category ${\rm Alg}(X)$
%${\rm Alg}(X_{[1]})$
of its algebra objects,
we obtain a functor
\[ {\rm Alg}: \laxmoncat{}\longrightarrow
              \cat \]
of $\infty$-categories.

Next, we consider coalgebra objects
of a monoidal $\infty$-category.
Let $p: X\to\Delta$ be a monoidal $\infty$-category,
that is,
$p$ is a Cartesian fibration 
such that the Segal morphisms are equivalences.
%For simplicity,
%we say that the underlying $\infty$-category
%$X_{[1]}$ is a monoidal $\infty$-category.
A coalgebra object of the monoidal $\infty$-category
$X_{[1]}$ is a section $s$ of $p$
such that $s$ takes inert morphisms of $\Delta$
to $p$-Cartesian morphisms.
As in the case of algebra objects,
we can regard a coalgebra object of $X_{[1]}$
as a morphism from ${\rm id}_{\Delta}$ to $p$
in $\oplaxmoncat{}$,
where ${\rm id}_{\Delta}: \Delta\to\Delta$
is a final object of $\oplaxmoncat{}$.

Let ${\rm Fun}_{\Delta}(\Delta,X)$
be the $\infty$-category of sections of $p$.
We denote by
\[ \coalg{}(X) \]
%\[ \coalg{}(X_{[1]}) \]
the full subcategory of 
${\rm Fun}_{\Delta}(\Delta,X)$
spanned by coalgebra objects of $X_{[1]}$.
Let $f$ be an oplax monoidal functor from
$p: X\to\Delta$ to $q: Y\to\Delta$.
Since $f$ carries $p$-Cartesian morphisms 
over inert morphisms of $\Delta$
to $q$-Cartesian morphisms,
it induces a functor
\[ \coalg{}(f): 
   \coalg{}(X)\longrightarrow \coalg{}(Y) \]
%\[ \coalg{}(F_{[1]}): 
%   \coalg{}(X_{[1]})\longrightarrow \coalg{}(Y_{[1]}) \]
of $\infty$-categories.
Assigning 
to a monoidal $\infty$-category $p: X\to\Delta$   
the $\infty$-category $\coalg{}(X)$
%$\coalg{}(X_{[1]})$
of its coalgebra objects,
we obtain a functor
\[ \coalg{}: \oplaxmoncat{}\longrightarrow
              \cat \]
of $\infty$-categories.
Note that the functor $\coalg{}$
is equivalent to
the composite of the functors
\[ \oplaxmoncat{}\stackrel{(-)^{\rm op}}{\longrightarrow}
   \laxmoncat{}\stackrel{{\rm Alg}}{\longrightarrow}
   \cat\stackrel{(-)^{\rm op}}{\longrightarrow}\cat.\]

\subsection{Bimonoids
in duoidal $\infty$-categories}
\label{subsection:bimonoida-in-duoidal-infty-category}
      
In this subsection
we introduce a notion of bimonoids
in duoidal $\infty$-categories.
For a duoidal $\infty$-category
$p: X\to \Delta_{\boxtimes}^{\rm op}\times\Delta_{\otimes}$,
we construct an $\infty$-category
$\bialg{}(X)$ of bimonoids in $X$.
This construction
determines a functor
$\bialg{}:\duo^{\rm bilax}\to\cat$
of $\infty$-categories.
Furthermore,
we show that the functor
$\coalg{}: \oplaxmoncat\to\cat$
can be upgraded to a functor
$\coalg{}{}^{\otimes}:\duo^{\rm bilax}\to \laxmoncat$,
and that 
the functor $\bialg{}$ is equivalent to the composite
$\alg{}\circ\coalg{}{}^{\otimes}$
of functors.

Recall that $\Delta_{\boxtimes}=\Delta_{\otimes}=\Delta$,
and $\pi_{\boxtimes}: 
\Delta_{\boxtimes}^{\rm op}\times\Delta_{\otimes}\to
\Delta_{\boxtimes}^{\rm op}$
and $\pi_{\otimes}:\Delta_{\boxtimes}^{\rm }\times\Delta_{\otimes}
\to\Delta_{\otimes}$
are the projections.
Note that the identity map 
${\rm id}_{\Delta_{\boxtimes}^{\rm op}\times\Delta_{\otimes}}: 
\Delta_{\boxtimes}^{\rm op}\times\Delta_{\otimes}\to
\Delta_{\boxtimes}^{\rm op}\times\Delta_{\otimes}$
is a final object of $\duo^{\rm bilax}$.
%Let $p: X\to\Delta_{\boxtimes}^{\rm op}\times \Delta_{\otimes}$
%be a duoidal $\infty$-category.
%Recall that $X_{[1],[1]}$
%is the underlying $\infty$-category of 
%the duoidal $\infty$-category $p$.
%For simplicity, we say that 
%the underlying $\infty$-category
%$X_{[1],[1]}$ is a duoidal $\infty$-cateogry.

\begin{definition}\rm
Let $p: X\to\Delta_{\boxtimes}^{\rm op}\times \Delta_{\otimes}$
be a duoidal $\infty$-category.
We define a bimonoid in the duoidal $\infty$-category
$X_{[1],[1]}$ to be
a morphism $B: {\rm id}_{\Delta_{\boxtimes}^{\rm op}
\times\Delta_{\otimes}}\to p$ in $\duo^{\rm bilax}$.
In other words,
a bimonoid in $X_{[1],[1]}$ is a
section $B$ of $p: X\to \Delta_{\boxtimes}^{\rm op}
\times\Delta_{\otimes}$:
\[ \xymatrix{
     & X\ar[d]^{p}\\
   \Delta_{\boxtimes}^{\rm op}\times \Delta_{\otimes}
   \ar[ur]^B\ar[r]^{\rm id}&
   \Delta_{\boxtimes}^{\rm op}\times \Delta_{\otimes}
}\]
such that $B$ sends 
$\pi_{\boxtimes}$-coCartesian 
%$(\pi_{\boxtimes}\circ {\rm Id})$-coCartesian 
morphisms over inert morphisms of $\Delta_{\boxtimes}^{\rm op}$
to $(\pi_{\boxtimes}\circ p)$-coCartesian morphisms, 
and that
$\pi_{\otimes}$-Cartesian morphisms 
%$(\pi_{\otimes}\circ {\rm Id})$-Cartesian morphisms 
over inert morphisms of $\Delta_{\otimes}$
to $(\pi_{\otimes}\circ p)$-Cartesian morphisms.

Let ${\rm Fun}_{\Delta_{\boxtimes}^{\rm op}\times\Delta_{\otimes}}
(\Delta_{\boxtimes}^{\rm op}\times\Delta_{\otimes},X)$
be the $\infty$-category of sections of 
$p:X\to\Delta_{\boxtimes}^{\rm op}\times\Delta_{\otimes}$.
We define an $\infty$-category
\[ \bialg{}(X) \]
to be the full subcategory
of ${\rm Fun}_{\Delta_{\boxtimes}^{\rm op}\times\Delta_{\otimes}}
(\Delta_{\boxtimes}^{\rm op}\times\Delta_{\otimes},X)$
spanned by bimonoids in $X_{[1],[1]}$,
and call it the $\infty$-category
of bimonoids of $X_{[1],[1]}$.
Assigning to a duoidal $\infty$-category 
$p: X\to\Delta_{\boxtimes}^{\rm op}\times\Delta_{\otimes}$
the $\infty$-category of bimonoids of $X_{[1],[1]}$,
we obtain a functor
\[ \bialg{}: \duo^{\rm bilax}\longrightarrow \cat \]
of $\infty$-categories.
\end{definition}

In the following of this subsection
we shall define a functor
\[ \coalg{}^{\otimes}: \duo^{\rm bilax}\longrightarrow \laxmoncat{}\]
and show that 
the functor
$\bialg{}: \duo^{\rm bilax}\longrightarrow\cat$
is equivalent to the composite of the functors
\[ \duo^{\rm bilax}\stackrel{\coalg{}^{\otimes}}
   {\hbox to 10mm{\rightarrowfill}}
   \laxmoncat{}\stackrel{\alg}
   {\hbox to 10mm{\rightarrowfill}}
   \cat.\]
By duality, 
we shall define a functor
\[ \alg^{\boxtimes}: \duo^{\rm bilax}\longrightarrow \oplaxmoncat{}\]
and show that
$\bialg{}$ is equivalent to
the composite of the functors
\[ \duo^{\rm bilax}\stackrel{\alg{}^{\boxtimes}}
    {\hbox to 10mm{\rightarrowfill}}
    \oplaxmoncat{}
    \stackrel{\coalg{}}
    {\hbox to 10mm{\rightarrowfill}}
    \cat.\]

%In the following of this section
%we take models of duoidal $\infty$-categories.
%Let $p: X\to\Delta_1^{\rm op}\times\Delta_{\otimes}$
%be a categorical fibration of simplicial sets,
%which represents a duoidal $\infty$-category in $\cat$.
Let $p: X\to\Delta_{\boxtimes}^{\rm op}\times\Delta_{\otimes}$
be a duoidal $\infty$-category.
We let $\overline{\coalg{}}{}^{\otimes}(X)$ 
be a simplicial set
over $\Delta_{\boxtimes}^{\rm op}$
satisfying the following formula
\[ {\rm Hom}_{{\rm Set}_{\Delta}/\Delta_{\boxtimes}^{\rm op}}
    (K,\overline{\coalg{}}{}^{\otimes}(X))\cong
   {\rm Hom}_{{\rm Set}_{\Delta}/(\Delta_{\boxtimes}^{\rm op}\times \Delta_{\otimes})}
    (K\times \Delta_{\otimes}, X) \]
for any simplicial set $K$ over $\Delta_{\boxtimes}^{\rm op}$,
where ${\rm Set}_{\Delta}$ is the category of simplicial sets.
We denote by $\overline{u}$
the map $\overline{\coalg{}}{}^{\otimes}(X)\to 
\Delta_{\boxtimes}^{\rm op}$.
Note that
the fiber $\overline{\coalg{}}{}^{\otimes}(X)_{[m]}$ of $\overline{u}$
at $[m]\in\Delta_{\boxtimes}^{\rm op}$
is ${\rm Fun}_{\Delta_{\otimes}}(\Delta_{\otimes},X_{[m],\bullet})$.

%\if0
\begin{lemma}
The map $\overline{u}:\overline{\coalg{}}{}^{\otimes}(X)
\to\Delta_{\boxtimes}^{\rm op}$
is a categorical fibration of simplicial sets.
In particular, $\overline{\coalg{}}{}^{\otimes}(X)$ 
is an $\infty$-category.
%is a quasi-category.
\end{lemma}

\proof
Since equivalences in $\Delta_{\boxtimes}^{\rm op}$ are 
only identities,
it suffices to show that 
$\overline{u}:\overline{\coalg{}}{}^{\otimes}(X)
\to\Delta_{\boxtimes}^{\rm op}$
is an inner fibration
by \cite[Corollary~2.4.6.5]{Lurie1}.
%[Lurie, HTT, Corollary~2.4.6.5]. 
Suppose we have a commutative diagram
\[ \xymatrix{
    \Lambda^n_i\ar[r]\ar[d] & 
    \overline{\coalg{}}{}^{\otimes}(X) \ar[d]\\%[2mm]
    \Delta^n \ar[r]\ar@{..>}[ur] & \Delta_{\boxtimes}^{\rm op}\\%[2mm]
   }\]
of simplicial sets for $0<i<n$.
We would like to construct a map
$\Delta^n\to \overline{\coalg{}}{}^{\otimes}(X)$,
which makes the whole diagram commutative.

By the definition of $\overline{\coalg{}}{}^{\otimes}(X)$,
the above diagram
is equivalent to the following commutative diagram
\[ \xymatrix{
    \Lambda^n_i\times\Delta_{\otimes} \ar[r]\ar[d]& 
    X \ar[d]\\
    \Delta^n\times\Delta_{\otimes} \ar[r]\ar@{..>}[ur]&
    \Delta_{\boxtimes}^{\rm op}\times\Delta_{\otimes}.\\
   }\]
The left vertical arrow is an inner anodyne map
by \cite[Corollary~2.3.2.4]{Lurie1}
%[Lurie, HTT, Corollary~2.3.2.4]
and the right vertical arrow is an inner fibration
by the assumption.
Hence there is a map $\Delta^n\times\Delta_{\otimes}\to X$
which makes the whole diagram commutative.
This induces a desired map
$\Delta^n\to \overline{\coalg{}}{}^{\otimes}(X)$.
\qed
%\fi

\begin{lemma}\label{lemma:in-order-to-show-cocartesian-fibration}
Suppose that we have a commutative diagram 
\[ \xymatrix{
    A \ar[rr]^f\ar[dr]_h && B\ar[dl]^g\\
     & C & \\
   }\]
of $\infty$-categories,
%of quasi-categories,
where $g$ and $h$ are coCartesian fibrations,
and 
$f$ is a categorical fibration
which carries $h$-coCartesian morphisms to $g$-coCartesian 
morphisms. 
Furthermore, suppose that we have a pullback diagram 
\[ \begin{array}{ccc}
    X & \stackrel{l}{\longrightarrow} & A\\[1mm]
    \mbox{\rm$\scriptstyle  k$}\bigg\downarrow
    \phantom{\mbox{\rm$\scriptstyle  k$}}
     & & 
    \phantom{\mbox{\rm$\scriptstyle  f$}}
    \bigg\downarrow\mbox{\rm$\scriptstyle  f$}\\[3mm]
    Y & \stackrel{m}{\longrightarrow} & B\\
   \end{array}\]
of $\infty$-categories.
%of quasi-categories.
If the image of any morphism of $Y$ under $m$
is a $g$-coCartesian morphism,
then $k: X\to Y$ is a coCartesian fibration. 
\end{lemma}

\proof
Note that $k$ is a categorical fibration
since it is a pullback of the categorical fibration
$f$ along $m$.
Let $x$ be an object of $X$
over $y\in Y$
and let $e: k(x)=y\to y'$
be a morphism of $Y$.
In order to prove that
$k$ is a coCartesian fibration,
we shall show that there is a
lifting $\overline{e}: x\to x'$ in $X$ of $e$, 
which is a $k$-coCartesian morphism. 

Since $h$ is a coCartesian fibration,
there is an $h$-coCartesian morphism
$\widetilde{e}: l(x)\to a$ over $g(m(e))$.
From the fact that $f$ carries $h$-coCartesian morphisms
to $g$-coCartesian morphisms,
$f(\widetilde{e}): f(l(x))=m(y)\to f(a)$
is a $g$-coCartesian morphism over $g(m(e))$.
By the assumption,
$m(e): m(y)\to m(y')$ is also a $g$-coCartesian morphism over $g(m(e))$,
and hence there is an equivalence
$e': f(a)\to m(y')$ 
in $B$
such that $e'\circ f(\widetilde{e})\simeq m(e)$.
By \cite[Corollary~2.4.6.5]{Lurie1},
we can take a lifting $\widetilde{e}': a\to a'$
in $A$ of $e'$,
which is an equivalence,
since $f$ is a categorical fibration
and $B$
is an $\infty$-category.
%is a quasi-category.
Using the fact that $f$ is an inner fibration,
we see that there is a $h$-coCartesian morphism
$\widehat{e}: l(x)\to a'$ in $A$,
which is a lifting of $m(e)$. 
By the dual of \cite[Proposition~2.4.1.3(3)]{Lurie1},
$\widehat{e}$ is an $f$-coCartesian morphism.
The lifting $\widehat{e}$ 
determines a lifting $\overline{e}: x\to x'$ of $e$,
which is a $k$-coCartesian morphism
by \cite[Proposition~2.4.1.3(2)]{Lurie1}.
This completes the proof.
\qed

\begin{lemma}\label{lemma:pbar-cocartesian-fibration}
The map $\overline{u}:
\overline{\coalg{}}{}^{\otimes}(X)\to\Delta_{\boxtimes}^{\rm op}$
is a coCartesian fibration of $\infty$-categories.
%In particular, $\overline{\coalg{}}(X)$ is a quasi-category.
\end{lemma}

\proof
We have a commutative triangle
\[ \xymatrix{
    {\rm Fun}(\Delta_{\otimes},X)\ar[rr]^{f}\ar[dr]_{h}&&
    {\rm Fun}(\Delta_{\otimes},\Delta_{\boxtimes}^{\rm op}\times\Delta_{\otimes})
    \ar[dl]^{g}\\
    & {\rm Fun}(\Delta_{\otimes},\Delta_{\boxtimes}^{\rm op})
   }\]
of $\infty$-categories, where
%of quasi-categories, where
$g$ and $h$ are coCartesian fibrations,
and the map $f$ carries $g$-coCartesian morphisms
to $h$-coCartesian morphisms
by \cite[Corollary~3.2.2.12]{Lurie1}.
We see that the map $f$ is a categorical fibration
by using \cite[Corollary~2.2.5.4]{Lurie1}.
Furthermore,
there is a pullback diagram
\[ \xymatrix{
    \overline{\coalg{}}{}^{\otimes}(X)
    \ar[r]\ar[d]_{\overline{u}} &
    {\rm Fun}(\Delta_{\otimes},X)\ar[d]^{f} \\
    \Delta_{\boxtimes}^{\rm op}\ar[r]^{j\hspace{10mm}} & 
    {\rm Fun}(\Delta_{\otimes},\Delta_{\boxtimes}^{\rm op}\times\Delta_{\otimes})\\
   } \]
of $\infty$-categories,
%of quasi-categories,
where the bottom horizontal arrow $j$ is 
an adjoint of the identity map of
$\Delta_{\boxtimes}^{\rm op}\times\Delta_{\otimes}$.
Since the image of any morphism of $\Delta_{\boxtimes}^{\rm op}$
under $j$
is a $g$-coCartesian morphism  
in ${\rm Fun}(\Delta_{\otimes},\Delta_{\boxtimes}^{\rm op}\times\Delta_{\otimes})$,
the lemma follows from
Lemma~\ref{lemma:in-order-to-show-cocartesian-fibration}.
%we see that the left vertical arrow is a coCartesian 
%fibration
%by using the dual of \cite[Proposition~2.4.1.3(3)]{Lurie1}.
\qed

\begin{remark}\rm
Let $e$ be a morphism of $\overline{\coalg{}}{}^{\otimes}(X)$,
which corresponds to a map $\Delta^1\times \Delta_{\otimes}\to X$.
The morphism $e$ is a $\overline{u}$-coCartesian 
morphism if and only if
the composite 
$\Delta^1\simeq \Delta^1\times\{[m]\}\hookrightarrow
\Delta^1\times\Delta_{\otimes}\to X$
is a $p$-coCartesian morphism
for all $[m]\in\Delta_{\otimes}$.
\end{remark}

\begin{definition}
We define 
\[ \coalg{}{}^{\otimes}(X) \]
to be
the full subcategory of $\overline{\coalg{}}{}^{\otimes}(X)$
spanned by coalgebra objects in $X_{[m],[1]}$
for some $[m]\in\Delta_{\boxtimes}^{\rm op}$.
By restriction of $\overline{u}$,
we obtain a map 
\[ u: \coalg{}{}^{\otimes}(X)
      \longrightarrow \Delta_{\boxtimes}^{\rm op} \]
of $\infty$-categories.
\end{definition}

We shall show that
$u: \coalg{}{}^{\otimes}(X)\to\Delta_{\boxtimes}^{\rm op}$
is a monoidal $\infty$-category.

\begin{proposition}\label{prop:cP-coCartesian-fibration}
The map $u: \coalg{}{}^{\otimes}(X)\to\Delta_{\boxtimes}^{\rm op}$
is a coCartesian fibration of $\infty$-categories.
\end{proposition}

\proof
Let $x$ be an object of $\coalg{}{}^{\otimes}(X)$ 
over $[m]\in\Delta_{\boxtimes}^{\rm op}$ and
let $e: [m]\to [m']$ be a morphism of $\Delta_{\boxtimes}^{\rm op}$.
Since $\overline{u}: \overline{\coalg{}}{}^{\otimes}(X)
\to\Delta_{\boxtimes}^{\rm op}$
is a coCartesian fibration by 
Lemma~\ref{lemma:pbar-cocartesian-fibration}, 
there is a $\overline{u}$-coCartesian morphism
$\overline{e}: x\to y$ over $e$.
In order to prove the proposition,
it suffices to show that
$y$ is an object of $\coalg{}{}^{\otimes}(X)$.
We have a map
$e_!:  X_{[m],\bullet}\to X_{[m'],\bullet}$
over $\Delta_{\otimes}$,
which carries $p_{[m],\bullet}$-Cartesian morphisms
over inner morphisms of $\Delta_{\otimes}$
to $p_{[m'],\bullet}$-Cartesian morphisms.
If we regard $x\in\coalg{}{}^{\otimes}(X)$
as a map $\Delta_{\otimes}\to X_{[m],\bullet}$,
then the object $y$ is identified with
the composite $e_!\circ x: \Delta_{\otimes}\to X_{[m'],\bullet}$
Since $x$ carries inert morphisms of $\Delta_{\otimes}$
to $p_{[m],\bullet}$-Cartesian morphisms,
$e_!\circ x$ carries inert morphisms of $\Delta_{\otimes}$
to $p_{[m'],\bullet}$-Cartesian morphisms.
Hence $y$ is an object of $\coalg{}{}^{\otimes}(X)$
over $[m']\in\Delta_{\boxtimes}^{\rm op}$.
This completes the proof.
\qed

\begin{theorem}\label{thm:cP-monoidal-infinity-category}
The map $u: \coalg{}{}^{\otimes}(X)\to\Delta_{\boxtimes}^{\rm op}$ 
is a monoidal $\infty$-category.
\end{theorem}

\proof
By Proposition~\ref{prop:cP-coCartesian-fibration},
it suffices to show that 
the Segal morphism
\begin{equation}\label{eq:Segal-morphism-coalg} 
   \coalg{}{}^{\otimes}(X)_{[m]}\longrightarrow
   \overbrace{\coalg{}{}^{\otimes}(X)_{[1]}
   \times\cdots\times
   \coalg{}{}^{\otimes}(X)_{[1]}}^m 
\end{equation}
is an equivalence of $\infty$-categories 
for each $[m]\in\Delta_{\boxtimes}^{\rm op}$.
Since $p: X\to\Delta_{\boxtimes}^{\rm op}\times \Delta_{\otimes}$
is a duoidal $\infty$-category,
the Segal morphism
$X_{[m],\bullet}\to
 X_{[1],\bullet}\times_{\Delta_{\otimes}}\cdots
              \times_{\Delta_{\otimes}}X_{[1],\bullet}$
%\overbrace{X_{[1],\bullet}\times_{\Delta_{\otimes}}\cdots
%              \times_{\Delta_{\otimes}}X_{[1],\bullet}}^m$
is an equivalence in $\cart{\Delta_{\otimes}}$
for each $[m]\in\Delta_{\boxtimes}^{\rm op}$.
This induces an equivalence
\[ \coalg{}{}^{\otimes}(X_{[m],\bullet})\longrightarrow
   \overbrace{\coalg{}{}^{\otimes}(X_{[1],\bullet})
   \times\cdots\times
   \coalg{}{}^{\otimes}(X_{[1],\bullet})}^m \]
%\[ \coalg{}(X_{[n],[1])\longrightarrow
%   \overbrace{\coalg{}(X_{[n],[1]})\times\cdots\times
%              \coalg{}(X_{[n],[1]})}^n \]
of $\infty$-categories.
Since we can identify this map with 
(\ref{eq:Segal-morphism-coalg}),
the map
$u:\coalg{}{}^{\otimes}(X)\to\Delta_{\boxtimes}^{\rm op}$
is a monoidal $\infty$-category. 
\qed

\bigskip

Next, we shall show that
a bilax monoidal functor
$f: X\to Y$ of duoidal $\infty$-categories
induces a lax monoidal functor
between the monoidal $\infty$-categories
$\coalg{}{}^{\otimes}(X)$
and $\coalg{}{}^{\otimes}(Y)$.
 
By the definition of
$\overline{\coalg{}}{}^{\otimes}(X)$,
we have an evaluation map
$\overline{\rm ev}:\overline{\coalg{}}{}^{\otimes}(X)
\times \Delta_{\otimes}\to X$
over $\Delta_{\boxtimes}^{\rm op}\times\Delta_{\otimes}$.
Let $f: X\to Y$ be a bilax monoidal functor between
duoidal $\infty$-categories.
Composing $\overline{\rm ev}$ with $f$,
we obtain a map 
$\overline{\coalg{}}{}^{\otimes}(X)\times \Delta_{\otimes}\to Y$
over $\Delta_{\boxtimes}^{\rm op}\times\Delta_{\otimes}$,
which induces a map 
$\overline{\coalg{}}{}^{\otimes}(f):
\overline{\coalg{}}{}^{\otimes}(X)\to
\overline{\coalg{}}{}^{\otimes}(Y)$
over $\Delta_{\boxtimes}^{\rm op}$.

Let $s$ be an object of 
$\coalg{}{}^{\otimes}(X)_{[m]}\simeq 
\coalg{}(X_{[m],\bullet})$
%$\coalg{}(X)_{[m]}\simeq \coalg{}(X_{[m],[1]})$
for $[m]\in\Delta_{\boxtimes}^{\rm op}$.
Since the bilax monoidal functor $f$ induces
an oplax monoidal functor
$f_{[m],\bullet}:X_{[m],\bullet}\to Y_{[m],\bullet}$,
we see that
$\overline{\coalg{}}{}^{\otimes}(f)(s)$
lands in the full subcategory
$\coalg{}{}^{\otimes}(Y)_{[m]}\simeq\coalg{}(Y_{[m],\bullet})$
%$\coalg{}(Y)_{[m]}\simeq\coalg{}(Y_{[m],[1]})$
of $\overline{\coalg{}}{}^{\otimes}(Y)_{[m]}$.
Thus, we obtain a functor
\[ \coalg{}{}^{\otimes}(f): 
   \coalg{}{}^{\otimes}(X)\longrightarrow 
   \coalg{}{}^{\otimes}(Y) \]
of $\infty$-categories over $\Delta_{\boxtimes}^{\rm op}$. 
 
\begin{proposition}\label{prop:bilax-induces-coalg-lax}
A bilax monoidal functor $f: X\to Y$ between
duoidal $\infty$-categories induces
a lax monoidal functor 
\[ \coalg{}{}^{\otimes}(f): 
   \coalg{}{}^{\otimes}(X)\longrightarrow
                \coalg{}{}^{\otimes}(Y)\]
between monoidal $\infty$-categories.
\end{proposition}

\proof
We have to show that
the induced map
$\coalg{}{}^{\otimes}(f):
\coalg{}{}^{\otimes}(X)\to\coalg{}{}^{\otimes}(Y)$
preserves coCartesian morphisms over
inert morphisms of $\Delta_{\boxtimes}^{\rm op}$.
Since $f: X\to Y$ is a bilax monoidal functor,
$f$ preserves coCartesian morphisms over
inert morphisms of $\Delta_{\boxtimes}^{\rm op}$.
If $e: [m]\to [m']$ is an inert morphism
of $\Delta_{\boxtimes}^{\rm op}$,
we have a commutative diagram
\[ \begin{array}{ccc}
    X_{[m],\bullet} & \stackrel{f_{[m],\bullet}}
    {\hbox to 15mm{\rightarrowfill}} & Y_{[m],\bullet}\\[2mm]
    \mbox{\rm $\scriptstyle e_!$}\bigg\downarrow
    \phantom{\mbox{\rm $\scriptstyle e_!$}}
     & & \phantom{\mbox{\rm $\scriptstyle e_!$}}
    \bigg\downarrow \mbox{\rm $\scriptstyle e_!$} \\[2mm]
    X_{[m'],\bullet} & \subrel{f_{[m'],\bullet}}
    {\hbox to 15mm{\rightarrowfill}} & Y_{[m'],\bullet} \\   
   \end{array}\]
of $\cart{(\Delta_{\otimes})^{\natural}}$.
This induces the following commutative diagram
\[ \begin{array}{ccc}
    \coalg{}(X_{[m],\bullet})& \stackrel{\coalg{}(f_{[m],\bullet})}
    {\hbox to 20mm{\rightarrowfill}} & \coalg{}(Y_{[m],\bullet})\\[2mm]
    \mbox{\rm $\scriptstyle e_!$}\bigg\downarrow
    \phantom{\mbox{\rm $\scriptstyle e_!$}}
     & & \phantom{\mbox{\rm $\scriptstyle e_!$}}
    \bigg\downarrow \mbox{\rm $\scriptstyle e_!$} \\[2mm]
    \coalg{}(X_{[m'],\bullet}) & \stackrel{\coalg{}(f_{[m'],\bullet})}
    {\hbox to 20mm{\rightarrowfill}} & \coalg{}(Y_{[m'],\bullet})\\[2mm]
   \end{array}\]
of $\infty$-categories.
This means that the functor
$\coalg{}{}^{\otimes}(f):
\coalg{}{}^{\otimes}(X)\to\coalg{}{}^{\otimes}(Y)$
preserves coCartesian morphisms
over inert morphisms of $\Delta_{\boxtimes}^{\rm op}$.
This completes the proof. 
\qed

\bigskip

By Theorem~\ref{thm:cP-monoidal-infinity-category} 
and 
Proposition~\ref{prop:bilax-induces-coalg-lax},
we obtain a functor
\[ \coalg{}{}^{\otimes}: \duo^{\rm bilax}\longrightarrow
             \laxmoncat{}\] 
by assigning to a duoidal $\infty$-category
$X\to\Delta_{\boxtimes}^{\rm op}\times\Delta_{\otimes}$
the monoidal $\infty$-category
$\coalg{}{}^{\otimes}(X)\to \Delta_{\boxtimes}^{\rm op}$.

\begin{remark}\rm
\label{remark:algoduo-bilax-oplax-monoidal-functor}
By duality,
there is a functor
\[ \alg{}^{\boxtimes}: \duo^{\rm bilax}\longrightarrow \oplaxmoncat{}\]
of $\infty$-categories.
Note that we have a commutative diagram
\[ \begin{array}{ccc}
     \duo^{\rm bilax} & 
     \stackrel{\coalg{}{}^{\otimes}}
     {\hbox to 10mm{\rightarrowfill}} & 
     \laxmoncat{} \\[2mm]
     \mbox{$\scriptstyle (-)^{\rm op}$}\bigg\downarrow
     \phantom{\mbox{$\scriptstyle (-)^{\rm op}$}}&&
     \phantom{\mbox{$\scriptstyle (-)^{\rm op}$}}\bigg\downarrow
     \mbox{$\scriptstyle (-)^{\rm op}$}\\[2mm]
     \duo^{\rm bilax} & 
     \stackrel{\alg{}^{\boxtimes}}
     {\hbox to 10mm{\rightarrowfill}} & 
     \oplaxmoncat{}\\
    \end{array}\]
in $\wcat$.
\end{remark}

By composing $\coalg{}{}^{\otimes}$ with
the functor $\alg{}: \laxmoncat{}\to\cat$,
we obtain a functor 
\[ \alg{}\circ\coalg{}{}^{\otimes}:
   \duo^{\rm bilax}\longrightarrow\cat.\]
Unwinding the definitions,
we obtain the following theorem.

\begin{theorem}
\label{theorem:bimonoid-is-algebra-of-coalgebra}
There is a natural isomorphism 
\[ \bialg{}(X)\simeq (\alg{}\circ\coalg{}{}^{\otimes})(X) \]
of simplicial sets
for any duoidal $\infty$-category
$p: X\to\Delta_{\boxtimes}^{\rm op}\times\Delta_{\otimes}$.
\end{theorem}

\begin{corollary}
\label{cor:equivalence-bimod-alg-coalg}
The functor $\bialg{}: \duo^{\rm bilax}\to\cat$
is equivalent to the composite $\alg{}\circ\coalg{}{}^{\otimes}$.
\end{corollary}

\begin{remark}\label{remark:coalg-alg-dual-duoidal}
\rm
By duality,
the functor $\bialg{}$ is equivalent to
the composite $\coalg{}\circ\alg{}^{\boxtimes}$.
Notice that there is a commutative diagram
\[ \begin{array}{ccccc}
     \duo^{\rm bilax} & \stackrel{\coalg{}{}^{\otimes}}
     {\hbox to 20mm{\rightarrowfill}} &
     \laxmoncat{} & \stackrel{\alg{}}
     {\hbox to 20mm{\rightarrowfill}} &
     \cat \\[2mm]
     \mbox{$\scriptstyle (-)^{\rm op}$}\bigg\downarrow
     \phantom{\mbox{$\scriptstyle (-)^{\rm op}$}}&&
     \mbox{$\scriptstyle (-)^{\rm op}$}\bigg\downarrow
     \phantom{\mbox{$\scriptstyle (-)^{\rm op}$}}&&
     \phantom{\mbox{$\scriptstyle (-)^{\rm op}$}}\bigg\downarrow
     \mbox{$\scriptstyle (-)^{\rm op}$}\\[2mm]
     \duo^{\rm bilax} & \stackrel{\alg{}^{\boxtimes}}
     {\hbox to 20mm{\rightarrowfill}} &
     \oplaxmoncat{} & \stackrel{\coalg{}}
     {\hbox to 20mm{\rightarrowfill}} &
     \cat \\[2mm]
   \end{array}\]
in $\wcat$.
\end{remark}

\subsection{Double monoids and double comonoids
in duoidal $\infty$-categories}
\label{subsection:double-monoid-comonoid-in-duoidal}

In this subsection
we introduce notions of double monoids
and double comonoids
in duoidal $\infty$-categories.
For a duoidal $\infty$-category
$p: X\to \Delta_{\otimes}^{\rm op}\times\Delta_{\boxtimes}^{\rm op}$,
we construct an $\infty$-category
$\dalg{}(X)$ of double monoids in $X$.
This construction
determines a functor
$\dalg{}:\duo^{\rm dlax}\to\cat$
of $\infty$-categories.
Furthermore, 
we show that the functor
$\alg{}: \laxmoncat\to\cat$
can be upgraded to a functor
$\alg{}{}^{\boxtimes}:\duo^{\rm dlax}\to \laxmoncat$,
and that 
the functor $\dalg{}$ is equivalent to the composite
$\alg{}\circ\alg{}{}^{\boxtimes}$.
By duality,
we also show the similar results on double comonoids.

First,
we define double monoids in a duoidal $\infty$-category.

\begin{definition}\rm
Let $p: X\to\Delta_{\otimes}^{\rm op}
\times \Delta_{\boxtimes}^{\rm op}$
be a duoidal $\infty$-category,
that is, it is a categorical fibration 
satisfying Condition~\ref{condition:duoidal-monoid-lax}.
We define a double monoid in the duoidal $\infty$-category
$X_{[1],[1]}$ to be
a morphism $D: {\rm id}_{\Delta_{\otimes}^{\rm op}
\times\Delta_{\boxtimes}^{\rm op}}\to p$ in $\duo^{\rm dlax}$.
In other words,
a double monoid in $X_{[1],[1]}$ is a
section $D$ of $p: X\to \Delta_{\otimes}^{\rm op}
\times\Delta_{\boxtimes}^{\rm op}$:
\[ \xymatrix{
     & X\ar[d]^{p}\\
   \Delta_{\otimes}^{\rm op}\times \Delta_{\boxtimes}^{\rm op}
   \ar[ur]^D\ar[r]^{\rm id}&
   \Delta_{\otimes}^{\rm op}\times \Delta_{\boxtimes}^{\rm op}
}\]
such that $D$ sends 
$\pi_{\otimes}$-coCartesian 
morphisms over inert morphisms of $\Delta_{\otimes}^{\rm op}$
to $(\pi_{\otimes}\circ p)$-coCartesian morphisms, 
and that 
$\Delta_{\boxtimes}^{\rm op}\simeq
\{[n]\}\times \Delta_{\boxtimes}^{\rm op}\to X_{[n],\bullet}$
is an algebra object of $X_{[n],[1]}$
for each $[n]\in\Delta_{\otimes}^{\rm op}$. 

Let ${\rm Fun}_{\Delta_{\otimes}^{\rm op}\times\Delta_{\boxtimes}^{\rm op}}
(\Delta_{\otimes}^{\rm op}\times\Delta_{\boxtimes}^{\rm op},X)$
be the $\infty$-category of sections of 
$p:X\to\Delta_{\otimes}^{\rm op}\times\Delta_{\boxtimes}^{\rm op}$.
We define an $\infty$-category
\[ \dalg{}(X) \]
to be the full subcategory
of ${\rm Fun}_{\Delta_{\otimes}^{\rm op}\times\Delta_{\boxtimes}^{\rm op}}
(\Delta_{\otimes}^{\rm op}\times\Delta_{\boxtimes}^{\rm op},X)$
spanned by double monoids in $X_{[1],[1]}$,
and call it the $\infty$-category
of double monoids of $X_{[1],[1]}$.
Assigning to a duoidal $\infty$-category 
$p: X\to\Delta_{\otimes}^{\rm op}\times\Delta_{\boxtimes}^{\rm op}$
the $\infty$-category of double monoids of $X_{[1],[1]}$,
we obtain a functor
\[ \dalg{}: \duo^{\rm dlax}\longrightarrow \cat \]
of $\infty$-categories.
\end{definition}

Next,
in the same way as the construction
of the functor 
$\coalg{}{}^{\otimes}: \duo^{\rm bilax}\to
\laxmoncat$ 
in \S\ref{subsection:bimonoida-in-duoidal-infty-category},
we construct a functor
$\alg{}{}^{\boxtimes}: \duo^{\rm dlax}\to\laxmoncat$,
which is a lifting
of the functor
$\alg{}:\laxmoncat\to\cat$.

For a duoidal $\infty$-category
$p: X\to \Delta_{\otimes}^{\rm op}\times\Delta_{\boxtimes}^{\rm op}$,
we let $\overline{\alg{}}{}^{\boxtimes}(X)$ 
be a simplicial set over 
$\Delta_{\otimes}^{\rm op}$
satisfying the following formula
\[ {\rm Hom}_{{\rm Set}_{\Delta}/\Delta_{\otimes}^{\rm op}}
    (K,\overline{\alg{}}{}^{\boxtimes}(X))\cong
   {\rm Hom}_{{\rm Set}_{\Delta}/
    (\Delta_{\otimes}^{\rm op}\times \Delta_{\boxtimes}^{\rm op})}
    (K\times \Delta_{\boxtimes}^{\rm op}, X) \]
for any simplicial set $K$ over $\Delta_{\otimes}^{\rm op}$.
We denote by $\overline{u}$
the map $\overline{\alg{}}{}^{\boxtimes}(X)\to 
\Delta_{\otimes}^{\rm op}$.
Note that
the fiber $\overline{\alg{}}{}^{\boxtimes}(X)_{[n]}$ of $\overline{u}$
at $[n]\in\Delta_{\otimes}^{\rm op}$
is ${\rm Fun}_{\Delta_{\boxtimes}^{\rm op}}
(\Delta_{\boxtimes}^{\rm op},X_{[n],\bullet})$.
We define 
\[ \alg{}{}^{\boxtimes}(X) \]
to be the full subcategory of
$\overline{\alg{}}{}^{\boxtimes}(X)$
spanned by algebra objects of 
the monoidal $\infty$-category
$X_{[n],\bullet}\to\Delta_{\boxtimes}^{\rm op}$
for some $[n]\in\Delta_{\otimes}^{\rm op}$.
We define a map
\[ u: \alg{}{}^{\boxtimes}(X)\to\Delta_{\otimes}^{\rm op} \]
to be the restriction of $\overline{u}$.

By the same argument as in 
Theorem~\ref{thm:cP-monoidal-infinity-category}
and
Proposition~\ref{prop:bilax-induces-coalg-lax},
we obtain the following theorem.

\begin{theorem}\label{theorem:duoidal-alg-monoidal-double-formulation}
The map $u: \alg{}{}^{\boxtimes}(X)\to\Delta_{\otimes}^{\rm op}$
is a monoidal $\infty$-category.
A double lax monoidal functor
$f: X\to Y$ of duoidal $\infty$-categories
induces a lax
monoidal functor
\[ \alg{}{}^{\boxtimes}(f): \alg{}{}^{\boxtimes}(X)
    \longrightarrow \alg{}{}^{\boxtimes}(Y).\]
\end{theorem}

\begin{remark}\rm
By Remark~\ref{remark:coalg-alg-dual-duoidal},
we have a monoidal $\infty$-category
$\alg{}{}^{\boxtimes}(Y)\to \Delta_{\otimes}$,
which is a Cartesian fibration,
for a duoidal $\infty$-category
$q: Y\to \Delta_{\boxtimes}^{\rm op}\times\Delta_{\otimes}$.
Let $p: X\to \Delta_{\otimes}^{\rm op}\times
\Delta_{\boxtimes}^{\rm op}$
be an equivalent duoidal $\infty$-category to $q$.
By Theorem~\ref{theorem:duoidal-alg-monoidal-double-formulation},
we obtain a monoidal $\infty$-category
$\alg{}{}^{\boxtimes}(X)\to\Delta_{\otimes}^{\rm op}$,
which is a coCartesian fibration.
These two monoidal $\infty$-category are equivalent.
\end{remark}

By Theorem~\ref{theorem:duoidal-alg-monoidal-double-formulation},
we obtain a functor
\[ \alg{}{}^{\boxtimes}: \duo^{\rm dlax}\longrightarrow
             \laxmoncat{}.\] 
%by assigning to a duoidal $\infty$-category
%$X\to\Delta_{\boxtimes}^{\rm op}\times\Delta_{\otimes}^{\rm op}$
%the monoidal $\infty$-category
%$\alg{}{}^{\boxtimes}(X)\to \Delta_{\otimes}^{\rm op}$.
Unwinding the definitions,
we obtain the following theorem.

\begin{theorem}
\label{theorem:dalg-equivalent-alg-alg-boxtimes}
There is a natural isomorphism 
\[ \dalg{}(X)\simeq (\alg{}\circ\alg{}{}^{\boxtimes})(X) \]
of simplicial sets
for any duoidal $\infty$-category
$p: X\to\Delta_{\otimes}^{\rm op}\times\Delta_{\boxtimes}^{\rm op}$.
The functor $\dalg{}: \duo^{\rm dlax}\to\cat$
is equivalent to the composite 
%$\alg{}\circ\alg{}{}^{\boxtimes}$:
\[ \duo^{\rm dlax}
   \stackrel{\alg{}^{\boxtimes}}
   {\hbox to 10mm{\rightarrowfill}}
   \laxmoncat 
   \stackrel{\alg{}}
   {\hbox to 10mm{\rightarrowfill}}
   \cat.\]
\end{theorem}

By duality,
we define double comonoids in duoidal $\infty$-categories.
  
\begin{definition}\rm
Let $p: X\to\Delta_{\otimes}\times \Delta_{\boxtimes}$
be a duoidal $\infty$-category,
that is, it is a categorical fibration 
which satisfies the conditions
in Remark~\ref{remark:dual-Cartesian-description-duoidal-category}.
%such that $R(p): X^{\rm op}\to
%\Delta_{\otimes}^{\rm op}\times \Delta_{\boxtimes}^{\rm op}$
%satisfies Condition~\ref{condition:duoidal-monoid-lax}.
We define a double comonoid in the duoidal $\infty$-category
$X_{[1],[1]}$ to be
a morphism $C: {\rm id}_{\Delta_{\otimes}
\times\Delta_{\boxtimes}}\to p$ in $\duo^{\rm doplax}$.
In other words,
a double comonoid in $X_{[1],[1]}$ is a
section $C$ of $p: X\to \Delta_{\otimes}\times\Delta_{\boxtimes}$:
\[ \xymatrix{
     & X\ar[d]^{p}\\
   \Delta_{\otimes}\times \Delta_{\boxtimes}
   \ar[ur]^C\ar[r]^{\rm id}&
   \Delta_{\otimes}\times \Delta_{\boxtimes}
}\]
such that $C$ sends 
$\pi_{\otimes}$-Cartesian 
morphisms over inert morphisms of $\Delta_{\otimes}$
to $(\pi_{\otimes}\circ p)$-Cartesian morphisms, 
and that 
$\Delta_{\boxtimes}\simeq
\{[n]\}\times \Delta_{\boxtimes}\to X_{[n],\bullet}$
is a coalgebra object of $X_{[n],[1]}$
for each $[n]\in\Delta_{\otimes}$. 

Let ${\rm Fun}_{\Delta_{\otimes}\times\Delta_{\boxtimes}}
(\Delta_{\otimes}\times\Delta_{\boxtimes},X)$
be the $\infty$-category of sections of 
$p:X\to\Delta_{\otimes}\times\Delta_{\boxtimes}$.
We define an $\infty$-category
\[ \dcalg{}(X) \]
to be the full subcategory
of ${\rm Fun}_{\Delta_{\otimes}\times\Delta_{\boxtimes}}
(\Delta_{\otimes}\times\Delta_{\boxtimes},X)$
spanned by double comonoids in $X_{[1],[1]}$,
and call it the $\infty$-category
of double comonoids of $X_{[1],[1]}$.
Assigning to a duoidal $\infty$-category 
$p: X\to\Delta_{\otimes}\times\Delta_{\boxtimes}$
the $\infty$-category of double comonoids of $X_{[1],[1]}$,
we obtain a functor
\[ \dcalg{}: \duo^{\rm doplax}\longrightarrow \cat \]
of $\infty$-categories.
\end{definition}

We define a functor
\[ \coalg{}^{\boxtimes}:\duo^{\rm doplax}\longrightarrow\cat\]
by making the following diagram commute
\[ \begin{array}{ccc}
     \duo^{\rm doplax}&
     \stackrel{\coalg{}{}^{\boxtimes}}
     {\hbox to 15mm{\rightarrowfill}}&
     \oplaxmoncat\\[2mm]
     \mbox{$\scriptstyle (-)^{\rm op}$}
     \bigg\downarrow
     \phantom{\mbox{$\scriptstyle (-)^{\rm op}$}} 
     & & 
     \phantom{\mbox{$\scriptstyle (-)^{\rm op}$}} 
     \bigg\downarrow
     \mbox{$\scriptstyle (-)^{\rm op}$}
     \\[2mm]
     \duo^{\rm dlax}&
     \stackrel{\alg{}{}^{\boxtimes}}
     {\hbox to 15mm{\rightarrowfill}}&
     \laxmoncat
   \end{array}\] 
in $\wcat$.

By the dual of 
Theorem~\ref{theorem:dalg-equivalent-alg-alg-boxtimes},
we obtain the following theorem.

\begin{theorem}
\label{theorem:dcoalg-equivalent-coalg-coalg-boxtimes}
There is a natural isomorphism 
\[ \dcalg{}(X)\simeq (\coalg{}\circ\coalg{}{}^{\boxtimes})(X) \]
of simplicial sets
for any duoidal $\infty$-category
$p: X\to\Delta_{\otimes}\times\Delta_{\boxtimes}$.
The functor $\dcalg{}: \duo^{\rm dlax}\to\cat$
is equivalent to the composite 
\[ \duo^{\rm doplax}
   \stackrel{\coalg{}^{\boxtimes}}
   {\hbox to 10mm{\rightarrowfill}}
   \oplaxmoncat 
   \stackrel{\coalg{}}
   {\hbox to 10mm{\rightarrowfill}}
   \cat.\]
\end{theorem}

%\input{monoidal-bicategories}
%\input{looping}
%\input{XA}
%\input{functoriality}
%\input{bicategory-to-duoidal}
%\input{MD}
%\input{proof-bicategory-to-duoidal}
%%%%%\input{cat-non-symmetric}
%%%\input{delta-product-duoidal}
%%%\input{lax-lax-formulation}
%%%\input{higher-monoidal}
%%%\input{o-p-duoidal-categories}
%%%\input{o-product-duoidal-categories}
%%%\input{enriched_duoidal}

%\newpage

%{\footnotesize
%\input{ref}

%}

\end{document}